\documentclass[12pt]{amsart} 
\usepackage[mathscr]{eucal}
\usepackage{graphics}
\input epsf 
\parskip=\smallskipamount \hoffset -1.5cm
\newtheorem{theorem}{Theorem}[section]

\newtheorem{corollary}[theorem]{Corollary}

\newtheorem{proposition}[theorem]{Proposition}

\newtheorem{problem}{Problem}

\numberwithin{equation}{section}

\newcommand{\boldC}{{\bf C}}
\newcommand{\boldD}{{\bf D}}

\newcommand{\boldR}{{\bf R}}
\newcommand{\boldS}{{\bf S}}

\newcommand{\eD}{{\mathfrak D}}

\newcommand{\eF}{{\mathfrak F}}
\newcommand{\eG}{{\mathfrak G}}
\newcommand{\eH}{{\mathfrak H}}
\newcommand{\eK}{{\mathfrak K}}
\newcommand{\eL}{{\mathfrak L}}
\newcommand{\eM}{{\mathfrak M}}
\newcommand{\eN}{{\mathfrak N}}
\newcommand{\eP}{{\mathfrak P}}

\newcommand{\Skappa}{\boldS_\kappa}

\newcommand{\Sone}{\boldS_1}

\newcommand{\Skappaprime}{\boldS_{\kappa'}}

\newcommand{\SkappaFG}{{\boldS_\kappa(\eF,\eG)}}
\newcommand{\SkappaprimeFG}{{\boldS_{\kappa'}(\eF,\eG)}}
\newcommand{\ip}[2]{\left<#1,#2\right>}
\newcommand{\sqminus}{{\rm{sq}_-}}

\newcommand{\sqplusminus}{{\rm{sq}_{\pm}}}
\newcommand{\oh}{{\mathcal O}}
\newcommand{\dom}{\hbox{\rm{dom}}}

\newcommand{\lfg}{\eL(\eF,\eG)}
\newcommand{\lf}{\eL(\eF)}

\newcommand{\lk}{\eL(\eK)}

\newcommand{\ran}{\hbox{\rm{ran\vphantom{d}}}}

\begin{document}

\phantom{x}  

\phantom{x}

\phantom{x}

\phantom{x}

\phantom{x}

\title
{A Note on Interpolation \\ [8pt]
  in the Generalized Schur Class} \author{D.  Alpay} \author{T.
  Constantinescu} \author{A. Dijksma} \author{J.  Rovnyak}
\address{Department of Mathematics \\
  Ben-Gurion University of the Negev \\
  P. O. Box 653 \\
  84105 Beer-Sheva, Israel} \email{\tt dany@math.bgu.ac.il}
\address{Programs in Mathematical Sciences \\
  University of Texas at Dallas \\
  Box 830688, Richardson, TX 75083-0688, U. S. A.}
\email{tiberiu@utdallas.edu} \address{
  Department of Mathematics \\
  University of Groningen \\
  P. O. Box 800 \\
  9700 AV Groningen, The Netherlands} \email{\tt dijksma@math.rug.nl}

\dedicatory{To Harry Dym: teacher, colleague and friend, in
  appreciation and with best wishes for his 60-th birthday.}

\address{
  Department of Mathematics \\
  University of Virginia \\
  Charlottesville, Virginia 22903-3199, U. S. A.}  \email{\tt
  rovnyak@Virginia.EDU} \thanks{A. Dijksma is grateful to
  Mr.~Harry~T.~Dozor for supporting his research through a Dozor
  Fellowship at the Ben-Gurion University of the Negev, Beer-Sheva,
  Israel.  J.~Rovnyak is supported by NSF Grant DMS-9801016.}

\begin{abstract}
  Realization theory for operator colligations on Pontryagin spaces is
  used to study interpolation and factorization in generalized Schur
  classes.  Several criteria are derived which imply that a given
  function is almost the restriction of a generalized Schur function.
  The role of realization theory in coefficient problems is also
  discussed; a solution of an indefinite Carath\'eodory-Fej\'er
  problem is obtained, as well as a result that relates the number of
  negative (positive) squares of the reproducing kernels associated
  with the canonical coisometric, isometric, and unitary realizations
  of a generalized Schur function to the number of negative (positive)
  eigenvalues of matrices derived from their Taylor coefficients.
\end{abstract}

\maketitle

\thispagestyle{empty}

\section{Introduction}

Let $\eF$ and $\eG$ be Pontryagin spaces having the same negative
index.  For any integer~$\kappa \ge 0$, the {\bf generalized Schur
  class} $\SkappaFG$ is the set of functions $S(z)$ with values in
$\lfg$ which are holomorphic on some subregion $\Omega$ of the open
unit disk $\boldD$ such that the kernel
\begin{equation}\label{E:ks}
  K_S(w,z) = \frac{1-S(z)S(w)^*}{1-z\bar w}
\end{equation}
has $\kappa$ negative squares.  In the scalar case, that is, when
$\eF=\eG=\boldC$ is the space of complex numbers in the Euclidean
metric, we write simply~$\Skappa$.  Terminology and notation used here
follow \cite{ADRSbook}, where an account of the realization theory for
the classes $\SkappaFG$ may be found (it should be noted that the
definition of $\SkappaFG$ in \cite{ADRSbook} requires that functions
are holomorphic at the origin, and we do not require this now).  For
example, we write $\sqplusminus K_S$ for the number of
positive/negative squares of the kernel \eqref{E:ks} and $\eH(S)$ for
the associated reproducing kernel Pontryagin space.  The generalized
Schur classes were extensively studied by Kre{\u\i}n and Langer
\cite{MR54:11103} when the coefficient spaces are Hilbert spaces,
which is the main case of interest here; in particular, generalized
Schur functions are meromorphic in $\boldD$.  The literature on
interpolation in such classes includes well-known works of Takagi
\cite{Takagi}, Adamjan, Arov, and Kre{\u\i}n \cite{MR45:7505},
Kre{\u\i}n and Langer \cite{MR57:1173}, Nudel$'$man
\cite{MR82f:30033}, Ball and Helton \cite{BH}, and others.

In interpolation theory, we consider kernels of the form \eqref{E:ks}
on $\Omega\times\Omega$, where $\Omega$ is any subset of $\boldD$ and
may be finite.  Nonnegativity of the kernel in this case implies that
$S(z)$ is the restriction of a classical Schur function.  Such a
conclusion cannot quite be drawn when \eqref{E:ks} has $\kappa$
negative squares: there are ``nonanalytic'' functions $S(z)$ such that
\eqref{E:ks} has a finite number of negative squares.  The following
result is a slightly more elaborate form of the standard example of
this phenomenon \cite[p.  82]{ADRSbook}.

\begin{proposition}\label{blaschke}
  Let $\Omega$ be a subset of $\boldD$ containing the point $w_0$ and
  at least one other point.  Let $S(z) = 1$ or $0$ according as
  $z=w_0$ or $z \in\Omega\setminus\{w_0\}$. Then $\sqminus K_S = 1$,
  and the following statements are equivalent:
  
  $(1)$ The function $S(z)$ is the restriction of a function in
  $\Sone$.
  
  $(2)$ The set $\Omega$ is a Blaschke sequence.
\end{proposition}

Recall that $\Omega$ is a Blaschke sequence if its points can be
arranged in a finite or infinite sequence $z_1, z_2,\dots$ with $\sum
(1-|z_n|^2) < \infty$. The zero set of a classical Schur function is a
Blaschke sequence \cite[p.~64]{RR1994}.  In the non-Blaschke case,
$S(z)$ is the restriction to $\Omega\setminus\{w_0\}$ of a function in
$S_{\kappa'}$, where $\kappa' = 0$, namely the function identically
zero.

\begin{proof}
  Let $H^2$ be the Hardy space on the unit disk, $H^2_\Omega$ the
  space of restrictions of functions in $H^2$ to $\Omega$.  We view
  $H^2_\Omega$ as a Hilbert space in the inner product such that the
  restriction mapping from $H^2$ onto $H^2_\Omega$ is a partial
  isometry.  The reproducing kernel for $H^2$ is the Szeg\H{o} kernel
  $k(w,z) = 1/(1-z\bar w)$, and the reproducing kernel for
  $H^2_\Omega$ is its restriction $k_\Omega(w,z)$ to
  $\Omega\times\Omega$.  The identity
\begin{equation}\label{E:example4}
K_S(w,z) = k_\Omega(w,z) - \frac{S(z)}{\sqrt{1-|w_0|^2}}\;
\frac{\overline{S(w)}}{\sqrt{1-|w_0|^2}}, \qquad w,z\in\Omega,
\end{equation}
shows that $\sqminus K_S \le 1$ \cite[Theorem 3.3]{ADRSsurvey}.  When
$z_1 = w_0$ and $z_2 \in \Omega \setminus \{w_0\}$, then
$$
\begin{pmatrix}
  K_S(z_i,z_j)
\end{pmatrix}_{i,j=1}^2
=
\begin{pmatrix}
  0 & \dfrac{1}{1 - z_2 \bar w_0} \\[10pt]
  \dfrac{1}{1 - w_0 \bar z_2} & \dfrac{1}{1 - |z_2|^2}
\end{pmatrix}.
$$
The determinant of this matrix is negative, and so the matrix has
one negative eigenvalue.  Hence $\sqminus K_S = 1$.

If (1) holds, then $\Omega\setminus\{w_0\}$ is contained in the zero
set of a nontrivial function in $\Sone$, and we obtain (2) because
such functions are of bounded type in $\boldD$.  Conversely, assume
(2) and consider the function
$$
S_1(z)=\left(\gamma \, \frac{z-\alpha}{1- z \bar
    \alpha}\right)^{-1}B(z), \quad \gamma= \left( \frac{w_0-\alpha}{1-
    w_0 \bar \alpha} \right)^{-1} B(w_0),
$$
where $B(z)$ is a Blaschke product having simple zeros at the
points of the Blaschke sequence $\Omega\setminus \{ w_0\}$ (so that
$0\neq |B(w_0)|<1$) and the number $\alpha$ is chosen in
$\boldD\setminus \Omega$ such that $|\gamma|=1$.  Then $S_1(z)$ is the
product of the inverse of a Blaschke factor and a classical Schur
function which does not vanish at the zero of the Blaschke factor, and
so it belongs to $\Sone$.  Evidently, $S(z)=S_1(z)$, $z \in \Omega$.
\end{proof}

Similar phenomena appear in Ball and Helton \cite{BH}.  Again consider
the scalar case.  Let $S(z)$ be defined on a subset $\Omega$ of
$\boldD$, and let $H^2_\Omega$ be as in the proof of
Proposition~\ref{blaschke}.  If $S(z)$ is a multiplier for
\smash{$H^2_\Omega$} (that is, \smash{$Sh\in H^2_{\Omega}$} for each
$h\in H^2_{\Omega}$) and \eqref{E:ks} has $\kappa$ negative squares,
there is a classical Schur function $S_0(z)$ and a Blaschke product
$B(z)$ of order~$\kappa$ such that $B(z)S(z) = S_0(z)$ on~$\Omega$.
Thus in this case, there exists a function $\hat S(z)$ in
$\boldS_{\kappa'}$ for some $\kappa'\le\kappa$ such that $S(z) = \hat
S(z)$ for all but at must $\kappa$ points of $\Omega$.  By
Proposition~\ref{blaschke}, it may occur that this can only be
satisfied with $\kappa'<\kappa$.

In Section~\ref{S:interpolation} we use realization theory to obtain
criteria which imply that a given function is the restriction of a
generalized Schur function, provided that certain exceptional points
are omitted.  These results are related to factorization theorems for
operator-valued functions of the Leech type \cite{ADRSlivsic}: that
is, we are given partially defined operator-valued functions $A(z)$
and $B(z)$, and it is required to find a generalized Schur function
$S(z)$ such that $B(z) = A(z)S(z)$.

In Section~\ref{S:coefficients} we discuss coefficient problems, in
which realization theory also plays a role in establishing analyticity
(see Theorem~\ref{T:CG0}).  Here some of the results are restricted to
the scalar case. Necessary conditions for the existence of solutions
are derived by considering the three kernels associated with a
generalized Schur function in its canonical coisometric, isometric,
and unitary realizations \cite{ADRSbook}.  First we show that these
conditions are, in fact, equivalent (Theorem~\ref{C:equivalence}).
Then we provide a complete solution to an indefinite form of the
Carath\'eodory-Fej\'er problem in the scalar case. Key to this result
is the equivalence of two matrix extension problems, one involving
lower triangular Toeplitz matrices and the other Hermitian Toeplitz
matrices.

\section{Interpolation and factorization}\label{S:interpolation}

Our approach is based on the use of characteristic functions of
partially isometric operator colligations, and the interpolation and
factorization criteria that we obtain are dictated by what is needed
to construct the colligations.
  
A scalar example gives an idea of the nature of the conditions.  Let
$S(z)$ be a complex-valued function defined on a nonempty subset
$\Omega$ of $\boldD$ such that the kernel $\eqref{E:ks}$ has $\kappa $
negative squares, and let $\eH(S;\Omega)$ be the associated
reproducing kernel Pontryagin space.  Let $e_w$ be the characteristic
function of a point $w$ of $\Omega$ ($e_w(z) = \delta_{wz}$ for all $z
\in \Omega$, where $\delta$ is the Kronecker symbol).

(1) {\it If $\Omega$ is not a Blaschke sequence, a necessary condition
  for the interpolation of $S(z)$ by a function $\hat S(z)$ in
  $\Skappa$ is that $e_{w}\notin\eH(S;\Omega)$ for every
  ${w}\in\Omega$.}

For if interpolation is possible, $\eH(S;\Omega)$ is the set of
restrictions of functions in $\eH(\hat S)$ to $\Omega$.  If
$e_{w}\in\eH(S;\Omega)$ for some ${w}\in\Omega$, then
$\Omega\setminus\{{w}\}$ is contained in the zero set of a nontrivial
function in $\eH(\hat S)$, and hence $\Omega$ is a Blaschke sequence
because functions in $\eH(\hat S)$ are of bounded type (for example,
see \cite[Theorem 4.2.3(4)]{ADRSbook}).

(2) {\it The same necessary condition does not necessarily hold if
  $\Omega$ is a Blaschke sequence.}

For example, suppose that $\kappa=0$ and $\Omega = \{{w}_1, {w}_2\}$
consists of two distinct points.  If $S(z)=z$ on $\Omega$, then
$K_S(w,z)=1$ identically on $\Omega\times\Omega$ and $\eH(S;\Omega)$
is a one-dimensional space consisting of constant functions on
$\Omega$; in this case, $e_{w}\notin\eH(S;\Omega)$ for all
${w}\in\Omega$.  But if $S(z)=z^2$, then $K_S(w,z)=1+z\bar w$ on
$\Omega\times\Omega$ and $\eH(S;\Omega)$ is two-dimensional; in this
case $e_{w}\in\eH(S;\Omega)$ for all ${w}\in\Omega$.  While the
condition is not always necessary for interpolation, it turns out that
such a condition can be sufficient.

The result below is stated in the form of a factorization problem and
thus has a possible systems interpretation.  We are given an ``input''
in the form of an operator-valued function $A(z)$ defined on some set
$\Omega$, which may be finite or infinite, and a target ``output''
function $B(z)$ on the same set.  It is required to find a transfer
function $S(z)$ for such a system.  The functions $A(z)$ and $B(z)$
themselves need not be holomorphic.

\begin{theorem}\label{rectangular}
  Let $\eF$, $\eG$, $\eK$ be Hilbert spaces, and let $\Omega$ be a
  subset of the unit disk containing the point $w_0$.  Let $A(z)$ and
  $B(z)$ be functions on $\Omega$ with values in $\eL(\eG,\eK)$ and
  $\eL(\eF,\eK)$.  Assume that the kernel
  \begin{equation}
    \label{ABkernel2}
    K(w,z) = \frac{A(z)A(w)^* - B(z)B(w)^*}{1 - \bar w z}
  \end{equation}
  has $\kappa$ negative squares on $\Omega\times\Omega$, and let
  $\eH_K$ be the associated reproducing kernel Pontryagin space.  Let
  $\eM$ be the subspace of $\eH_K\oplus \eG$ consisting of all
  elements $k(z)\oplus g$ such that
  $$A(w_0)g=0 \quad\text{and}\quad \frac{z-w_0}{\sqrt{1-|w_0|^2}} \,
  k(z)+ [A(z)-A(w_0)]g \equiv 0 \quad\text{on}\quad \Omega.$$
  Let
  $\eN$ be the subspace of $\eH_K\oplus \eF$ consisting of all
  elements $h(z)\oplus f$ such that
  $$
  \frac{1-z \bar w_0}{\sqrt{1-|w_0|^2}}\, h(z)+B(z)f\equiv 0
  \quad\text{on}\quad \Omega.
  $$
  Assume that $\eM$ and $\eN$ are Hilbert spaces in the inner
  products of the larger spaces.  Then there is a function $S(z) \in
  \SkappaprimeFG$ for some $\kappa'\le\kappa$ such that $ B(z) =
  A(z)S(z)$ for $z=w_0$ and for all but at most $\kappa$ points $z$ of
  $\Omega \setminus \{w_0\}$.  In this case, $\kappa' = \kappa$ if and
  only if the elements $h$ of $\eH(S)$ such that $A(z)h(z) \equiv 0$
  on $\Omega$ form a Hilbert subspace of $\eH(S)$.
\end{theorem}

The function $S(z)$ which is constructed in the proof is holomorphic
at $w_0$.  The subspaces $\eM$ and $\eN$ defined in the statement of
the theorem are automatically closed by the continuity of function
values in a reproducing kernel space \cite[Theorem 1.1.2]{ADRSbook}.

\begin{proof}  It is sufficient to prove the result when $0 \in
  \Omega$ and $w_0 = 0$.  For suppose that the result is known in this
  case, and consider the general situation.  Let $\varphi$ be the
  linear fractional mapping of $\boldD$ onto itself given by
  $\varphi(z) = (w_0 - z)/(1 - \bar w_0 z)$.  Thus $\varphi(w_0)=0$
  and $\varphi^{-1} = \varphi$.  Put $\Omega' = \varphi(\Omega)$,
  $w_0' = 0$, and
  \begin{align}
    A'(z) &= A(\varphi^{-1}(z)), \qquad z \in \Omega',   \notag \\
    B'(z) &= B(\varphi^{-1}(z)), \qquad z \in \Omega'.  \notag
  \end{align}
  Define $K'(w,z)$ on $\Omega'\times\Omega'$ by \eqref{ABkernel1}
  using $A'(z)$ and $B'(z)$ in place of $A(z)$ and $B(z)$.  A short
  calculation shows that
  $$
  K'(w,z) = \frac{1 - |w_0|^2}{(1 - \bar w_0 z)(1 - w_0 \bar w)}
  \,K(\varphi^{-1}(w),\varphi^{-1}(z)), \qquad w,z \in \Omega',
  $$
  and so $\sqminus K' = \kappa$; write $\eH_{K'}$ for the
  associated reproducing kernel Pontryagin space.  The preceding
  reproducing kernel identity may be used to show that the mapping
  $$
  V' \colon f(z) \to \frac{\sqrt{1 - |w_0|^2}}{1 - \bar w_0 z}
  \;f(\varphi^{-1}(z))
  $$
  acts as an isometry from $\eH_{K}$ onto $\eH_{K'}$.  Writing
  $\eM$ and $\eN$ for the subspaces defined in the theorem for the
  original functions $A(z)$ and $B(z)$ and point $w_0 \in \Omega$, and
  $\eM'$ and $\eN'$ for the corresponding subspaces relative to
  $A'(z)$ and $B'(z)$ and point $w_0' \in \Omega'$, we find that
  $$
  (V' \oplus -1_{\eG})\, \eM = \eM' \quad\text{and}\quad (V' \oplus
  1_{\eF})\, \eN = \eN.
  $$
  Since we assume the result when $0 \in \Omega$ and $w_0 = 0$, we
  can find a function $S'(z) \in \SkappaprimeFG$ for some
  $\kappa'\le\kappa$ such that $S'(z)$ is holomorphic at $w_0'=0$ and
  $B'(z) = A'(z)S'(z)$ for $z=w_0'$ and for all but at most $\kappa$
  points $z$ of $\Omega' \setminus \{w_0'\}$.  Then $S(z) =
  S'(\varphi(z))$ has the required properties.
  
  Thus without loss of generality, we may assume that $0 \in \Omega$
  and $w_0 = 0$.  Define a linear relation $\boldR$ in
  $(\eH_K\oplus\eG) \times (\eH_K\oplus\eF)$ as the span of all pairs
\begin{equation}\label{generatingpairs}
{\pmatrix 
\pmatrix 
\vphantom{\dfrac{K(\alpha,\cdot) - K(0,\cdot)}{\bar \alpha}} 
K(\alpha,\cdot)u_1 \\
\noalign{\vskip4pt} 
\vphantom{\dfrac{\overline{S(\alpha)} -
\overline{S(0)}}{\bar \alpha}}
\dfrac{A(\alpha)^* - A(0)^*}{\bar \alpha}
\, u_1 + A(0)^*u_2
\endpmatrix , 
&\!\!\!
\pmatrix
\dfrac{K(\alpha,\cdot) - K(0,\cdot)}{\bar \alpha}  \,
u_1 + K(0,\cdot) u_2 \\ \noalign{\vskip6pt}
\dfrac{B(\alpha)^* - B(0)^*}{\bar \alpha}
\, u_1 + B(0)^*u_2 
\endpmatrix 
\endpmatrix}
\end{equation}
with $\alpha\in\Omega \setminus \{0\}$ and $u_1,u_2\in\eK$.  A direct
calculation shows that $\boldR$ is isometric.  In fact, consider a
second pair with $\alpha$ replaced by $\beta$ and $u_1,u_2$ replaced
by $v_1,v_2$.  Expand and simplify the inner products of the first
members in $\eH_K\oplus\eG$ and second members in $\eH_K\oplus\eF$.
After simplification, in both cases we obtain
\begin{align}
  \label{simplified}
  {\ip{K(\alpha,\beta)u_1}{v_1}}_{\eK} &+
  {\ip{\frac{A(\beta)A(\alpha)^* - A(\beta)A(0)^* - A(0)A(\alpha)^*
        + A(0)A(0)^*}{\bar \alpha \beta}\,u_1}{v_1}}_{\eK} \\
  &\hskip-2cm + {\ip{\frac{A(0)A(\alpha)^* - A(0)A(0)^*}{\bar
        \alpha}\,u_1}{v_2}}_{\eK} + {\ip{\frac{A(\beta)A(0)^*
        - A(0)A(0)^*}{\beta}\,u_2}{v_1}}_{\eK}  \notag \\[6pt]
  &\hskip-2cm+ {\ip{A(0)A(0)^*u_2}{v_2}}_{\eK}, \notag
\end{align}
and this verifies the assertion.  The orthogonal complement of the
domain of $\boldR$ is $\eM$, and the orthogonal complement of the
range of $\boldR$ is $\eN$.  Since these are Hilbert spaces, it
follows from \cite[Theorem 1.4.2]{ADRSbook} that there is a continuous
partial isometry
$$
V = \left(\begin{matrix} T&F\\G&H \end{matrix}\right): \pmatrix
\eH_K\\ \eF \endpmatrix \rightarrow \pmatrix \eH_K \\ \eG \endpmatrix
$$
such that $V^*$ has initial space $\overline{\dom}\, \boldR$ and
final space $\overline{\ran}\, \boldR$ and
\begin{equation}\notag
V^*: \pmatrix 
\vphantom{\dfrac{K(\alpha,\cdot) -K(0,\cdot)}{\bar \alpha}} 
K(\alpha,\cdot)u_1 \\
\noalign{\vskip4pt} 
\vphantom{ \dfrac{ S(\alpha) - \overline{S(0)} }{\bar \alpha}}
\dfrac{A(\alpha)^* - A(0)^*}{\bar \alpha}
\, u_1 + A(0)^*u_2
\endpmatrix  
\rightarrow
\pmatrix
\dfrac{K(\alpha,\cdot) - K(0,\cdot)}{\bar \alpha} \,
u_1 + K(0,\cdot) u_2 \\ \noalign{\vskip6pt}
\dfrac{B(\alpha)^* - B(0)^*}{\bar \alpha}
\, u_1 + B(0)^*u_2 \,
\endpmatrix 
\end{equation}
for all $\alpha\in\Omega\setminus\{0\}$ and $u_1,u_2\in\eK$.
Calculating as in \cite[p.~51]{ADRSbook}, we find that
\begin{align}
  (Th)(z) &= \frac{h(z)-A(z)Gh}{z}, \hskip1.5cm z\in
  \Omega\setminus\{0\},
  \notag\\
  (Ff)(z) &= \frac{B(z)-A(z)H}{z}\,f, \hskip1.3cm z\in
  \Omega\setminus\{0\},
  \notag\\
  A(0)Gh &= h(0) , \notag\\
  A(0)Hf &= B(0)f , \notag
\end{align}
for all $h\in\eH_K$ and $f\in\eF$.

Since $V$ is a partial isometry whose kernel is a Hilbert space, $V$
is a contraction.  The embedding mappings $E_\eF$ and $E_\eG$ from
$\eH_K$ into $\eH_K\oplus \eF$ and $\eH_K\oplus \eG$ are contractions
(in fact isometries), and their adjoints act as projections.  The
adjoints are also contractions because we assume that $\eF$ and $\eG$
are Hilbert spaces.  Therefore
$$
T = E_\eG^* V E_\eF
$$
is a contraction on the Pontryagin space $\eH_K$.  By \cite[Lemma
11.1 (p.~75)]{IKL}, the part of the spectrum of $T$ that lies in
$|\lambda|>1$ consists of normal eigenvalues.  By \cite[Theorem 11.2
(p.~84)]{IKL}, the span of root manifolds for eigenvalues in
$|\lambda|>1$ is contained in a nonpositive subspace, and hence the
number of such eigenvalues is at most $\sqminus \eH_K = \kappa$.  It
follows that $1-zT$ is invertible for all but at most $\kappa$ points
in $\boldD$.  Since these exceptional points obviously do not include
$0$, $1-zT$ is invertible for all $z\in \Omega\setminus
\{\lambda_1,\dots,\lambda_q \}$ for some nonzero numbers
$\lambda_1,\dots,\lambda_q$ in $\boldD$; here $q\le\kappa$ and
possibly $q=0$ when there are no exceptional points.

\medskip
\noindent
{\it Claim $1$}: {\it If $w\in \Omega\setminus
  \{\lambda_1,\dots,\lambda_q \}$, $h\in\eH_K$, and $(1 -
  wT)^{-1}h=g$, then
\begin{equation}\label{eq1}
g(z) = \frac{zh(z)-wA(z)Gg}{z-w},
\qquad z\in\Omega\setminus\{w\},
\end{equation}
and $ h(w) = A(w)Gg$. }

Since this is trivially true if $w=0$, assume that $w\neq0$.  Then
\begin{equation}\label{eq2} 
h(z) = g(z) - w\,\frac{g(z)-A(z)Gg}{z},\qquad
z\in\Omega\setminus\{0\}.
\end{equation}
Since $w\neq 0$, we can take $z=w$ in \eqref{eq2} to get
$h(w)=A(w)Gg$.  Again by \eqref{eq2},
$$
(z-w)g(z) = zh(z)-wA(z)Gg
$$
for $z\in\Omega\setminus\{0\}$.  Trivially the last identity holds
for $z=0$ as well, and we obtain \eqref{eq1}.

\medskip
\noindent {\it Claim $2$}: {\it Define $S(w) = H + wG(1-wT)^{-1}F$
  for all $w \in \boldD\setminus \{\lambda_1,\dots,\lambda_q \}$.
  Then
  $$B(w) = A(w)S(w)$$
  for all $w \in \Omega\setminus
  \{\lambda_1,\dots,\lambda_q \}$.  }

The case $w=0$ is clear.  Assume $w \in \Omega\setminus
\{\lambda_1,\dots,\lambda_q \}$ and $w\neq0$. Fix $f\in\eF$.  We use
Claim~$1$ with $g = (1-wT)^{-1}h$, $h = Ff$.  Thus
\begin{equation}\notag
wA(w)G(1-wT)^{-1}Ff = wA(w)Gg = wh(w) = w\, (Ff)(w)
 = B(w)f - A(w)Hf.
\end{equation}
Claim~$2$ follows.

\goodbreak \medskip
\noindent {\it Claim $3$}: {\it $S\in\Skappaprime$ for some
  $\kappa'\le\kappa$.}

It is clear from the definition of $S(z)$ that it is a holomorphic
function on $\boldD\setminus \{\lambda_1,\dots,\lambda_q \}$.  For all
$w,z\in \boldD\setminus \{\lambda_1,\dots,\lambda_q \}$, by the
identity \cite[(1.2.9)]{ADRSbook}, {\allowdisplaybreaks
\begin{align}
  1 - S(z)S(w)^* &=
\begin{pmatrix}
  G(1 -zT)^{-1} &1
\end{pmatrix}
\begin{pmatrix}
  (1 -\bar w T^*)^{-1}G^* \notag \\
  1
\end{pmatrix}
\\
&\qquad - \begin{pmatrix} zG(1 -zT)^{-1} &1
\end{pmatrix}
VV^*
\begin{pmatrix}
  \bar w (1 -\bar w T^*)^{-1}G^* \\ 1
\end{pmatrix}  \notag \\
&=
\begin{pmatrix}
  G(1 -zT)^{-1} &1
\end{pmatrix}
\begin{pmatrix}
  (1 -\bar w T^*)^{-1}G^* \notag \\ 1
\end{pmatrix} \\
&\qquad - \begin{pmatrix} zG(1 -zT)^{-1} &1
\end{pmatrix}
\begin{pmatrix}
  \bar w (1 -\bar w T^*)^{-1}G^*
  \\
  1
\end{pmatrix}
\notag \\
&\qquad +
\begin{pmatrix}
  zG(1 -zT)^{-1} &1
\end{pmatrix}
(1-VV^*)
\begin{pmatrix}
  \bar w (1 -\bar w T^*)^{-1}G^*
  \\
  1
\end{pmatrix}
\notag \\
& = (1-z\bar w) G(1-zT)^{-1}(1-\bar w T^*)^{-1}G^*
\notag \\*
&\qquad + \begin{pmatrix} zG(1 -zT)^{-1} &1
\end{pmatrix}
(1-VV^*)
\begin{pmatrix}
  \bar w (1 -\bar w T^*)^{-1}G^*
  \\
  1
\end{pmatrix}
.  \notag
\end{align}
} $\!\!$Since $1-VV^* \ge 0$ in the partial ordering of selfadjoint
operators, $1-VV^* = MM^*$ for some operator
$M\in\eL(\eD,\eH_K\oplus\eG)$, where $\eD$ is a Hilbert space (see,
for example, \cite[Theorem 2.1]{DR1996}; we can choose $M$ so that it
has zero kernel, but this property is not needed.  Therefore
\begin{equation}\label{E:kshat}
K_{S}(w,z) = G(1-zT)^{-1}(1-\bar w T^*)^{-1}G^* + 
\frac{\Phi(z)\Phi(w)^*}{1-z\bar w},
\qquad w,z\in 
\boldD\setminus \{\lambda_1,\dots,\lambda_q \},
\end{equation}
where
$$
\Phi(z) = \pmatrix zG(1 -zT)^{-1} & 1 \endpmatrix M, \qquad z\in
\boldD \setminus \{\lambda_1,\dots,\lambda_q \},
$$
is a holomorphic function with values in $\eL(\eD,\eG)$.  The first
summand on the right of \eqref{E:kshat} has $\kappa''$ negative
squares for some $\kappa''\le\kappa$ by \cite[Lemma
$1.1.1'$]{ADRSbook}, and the second summand is nonnegative because
$\eD$ is a Hilbert space.  Thus by \cite[Theorem 1.5.5]{ADRSbook} the
kernel \eqref{E:kshat} has $\kappa'$ negative squares, where
$\kappa'\le\kappa''\le\kappa$.  Hence $S\in\Skappaprime$, which proves
Claim 3.

The function $S(z)$ has the required properties by Claims~2 and~3.
The last statement, which gives the condition for $\kappa' = \kappa$,
follows from \cite[Theorem 1.5.7]{ADRSbook}.
\end{proof}

The next result identifies a case in which the conditions in
Theorem~\ref{rectangular} can be verified.  Namely, we assume that the
values of $A(z)$ are ``square'' in the sense that $\eK=\eG$ and so the
values of $A(z)$ are in $\eL(\eG)$.  We also assume that one of these
values is invertible, and we take this to be $1_\eG$.

\goodbreak
\begin{theorem}\label{square}
  Let $\eF$ and $\eG$ be Hilbert spaces, and let $A(z)$ and $B(z)$ be
  functions which are defined on a subset $\Omega$ of $\boldD$ with
  values in $\eL(\eG)$ and $\lfg$.  Assume that the kernel
  \begin{equation}
    \label{ABkernel1}
    K(w,z) = \frac{A(z)A(w)^* - B(z)B(w)^*}{1 - \bar w z}
  \end{equation}
  has $\kappa$ negative squares on $\Omega\times\Omega$, and let
  $\eH_K$ be the associated reproducing kernel Pontryagin space.
  Assume that there is a point $w_0\in\Omega$ such that
  \begin{enumerate}
  \item[(1)] $A(w_0)=1_\eG$, and
  \item[(2)] the set of elements of $\eH_K$ which vanish on $\Omega
    \setminus \{w_0\}$ is a Hilbert subspace of $\eH_K$.
\end{enumerate}
Then there is a function $S(z) \in \SkappaprimeFG$ for some
$\kappa'\le\kappa$ such that $ B(z) = A(z)S(z)$ for $z=w_0$ and for
all but at most $\kappa$ points $z$ of $\Omega \setminus \{w_0\}$.  In
this case, $\kappa' = \kappa$ if and only if the elements $h$ of
$\eH(S)$ such that $A(z)h(z) \equiv 0$ on $\Omega$ form a Hilbert
subspace of $\eH(S)$.
\end{theorem}

The function $S(z)$ constructed in the proof is holomorphic at $w_0$.

\begin{proof}
  The last statement follows from \cite[Theorem 1.5.7]{ADRSbook}.  It
  is convenient to assume that $0 \in \Omega$ and $w_0 = 0$.  If the
  result is known in this case, then as in the proof of
  Theorem~\ref{rectangular}, define $A'(z)$ and $B'(z)$ on $\Omega' =
  \varphi(\Omega)$, where $\varphi(z) = (w_0-z)/(1 - \bar w_0 z)$.  As
  in the same proof, introduce the kernel $K'(w,z)$ and isomorphism
  $V'$ from $\eH_K$ onto $\eH_{K'}$.  Under $V'$, the functions in
  $\eH_{K}$ which vanish on $\Omega \setminus \{w_0\}$ correspond to
  the functions in $\eH_{K'}$ which vanish on $\Omega' \setminus
  \{w_0'\}$, where $w_0' = 0$.  Then as before, the special case
  implies the general result.
  
  In what follows, we assume that $0 \in \Omega$ and $w_0 = 0$.  We
  apply Theorem~\ref{rectangular} in this situation and also with $\eK
  = \eG$.  It is easy to see that the subspace $\eM$ in
  Theorem~\ref{rectangular} coincides with the the set of elements of
  $\eH_K$ which vanish on $\Omega \setminus \{0\}$ and is thus a
  Hilbert space by hypothesis.  We show that the subspace $\eN$ in
  Theorem~\ref{rectangular} is a Hilbert space.  By the first part of
  the proof of Theorem~\ref{rectangular}, $\eN$ is the orthogonal
  complement of the range of the relation $\boldR$ in $\eH_K \oplus
  \eF$, and therefore it is the same thing to show that the range of
  $\boldR$ contains a strictly negative subspace of
  dimension~$\kappa$.  By \cite[Lemma 1.1.1$'$]{ADRSbook}, it is
  sufficient to show that some Gram matrix of elements of the range of
  $\boldR$ has $\kappa$ negative eigenvalues.  In fact, consider two
  of the second members of the pairs \eqref{generatingpairs} that
  define $\boldR$, say
  $$
  \pmatrix \dfrac{K(\alpha,\cdot) - K(0,\cdot)}{\bar \alpha} \, u_1
  + K(0,\cdot) u_2 \\ \noalign{\vskip6pt} \dfrac{B(\alpha)^* -
    B(0)^*}{\bar \alpha} \, u_1 + B(0)^*u_2 \endpmatrix \quad
  \hbox{\rm and} \quad \pmatrix \dfrac{K(\beta,\cdot) -
    K(0,\cdot)}{\bar \beta} \, v_1 + K(0,\cdot) v_2 \\
  \noalign{\vskip6pt} \dfrac{B(\beta)^* - B(0)^*}{\bar \beta} \, v_1 +
  B(0)^*v_2 \endpmatrix .
  $$
  By \eqref{simplified}, since now $A(0) = 1_\eG$, the inner
  product of these elements in $\eH_K \oplus \eF$ is equal to
\begin{align}
  {\ip{K(\alpha,\beta)u_1}{v_1}}_{\eG} &+
  {\ip{\frac{A(\beta)A(\alpha)^* - A(\beta) - A(\alpha)^*
        + 1_\eG}{\bar \alpha \beta}\,u_1}{v_1}}_{\eG} \notag \\
  &\hskip-2cm + {\ip{\frac{A(\alpha)^* - 1_\eG}{\bar
        \alpha}\,u_1}{v_2}}_{\eG} + {\ip{\frac{A(\beta) -
        1_\eG}{\beta}\,u_2}{v_1}}_{\eG}
  + {\ip{u_2}{v_2}}_{\eG} \notag \\[8pt]
  &\hskip-2cm = {\ip{K(\alpha,\beta)u_1}{v_1}}_{\eG} +
  {\ip{\frac{A(\alpha)^* - 1_\eG}{\bar \alpha}\,u_1 + u_2}
    {\frac{A(\beta)^* - 1_\eG}{\bar \beta}\,v_1 + v_2}}_{\eG}.  \notag
\end{align}
Here we can choose $\alpha,\beta$ and $u_1,u_2$ arbitrarily, and then
choose $v_1,v_2$ so that
$$
\frac{A(\alpha)^* - 1_\eG}{\bar \alpha}\,u_1 + u_2 =
\frac{A(\beta)^* - 1_\eG}{\bar \beta}\,v_1 + v_2 = 0.
$$
Since we assume that $\sqminus K = \kappa$, it follows that some
Gram matrix of elements of the range of $\boldR$ has $\kappa$ negative
eigenvalues, as was to be shown.  This completes the proof that $\eN$
is a Hilbert space.

The hypotheses of Theorem~\ref{rectangular} are thus met, and
Theorem~\ref{rectangular} yields a function $S(z) \in \SkappaprimeFG$,
$\kappa'\le\kappa$, such that $ B(z) = A(z)S(z)$ for $z=0$ and for all
but at most $\kappa$ points $z$ of $\Omega \setminus \{0\}$.
\end{proof}

We give another condition for interpolation.  Suppose that $S(z)$
belongs to $\Skappa$ and is holomorphic at the origin.  Then $zS(z)$
also belongs to $\Skappa$, and thus both kernels
\begin{equation}\notag
\frac{1 - S(z)\overline{S(w)}}{1 - z\bar w}
\qquad \text{and}\qquad
\frac{1 - z\bar w S(z)\overline{S(w)}}{1 - z\bar w}
\end{equation}
have $\kappa$ negative squares (see \cite[Example~1 on
p.~132]{ADRSbook}).  In the other direction, a condition on two
kernels is sufficient for interpolation from an arbitrary set $\Omega$
with at most a finite number of exceptional points.

\goodbreak
\begin{theorem}\label{th3}
  Let $A(z)$ and $B(z)$ be functions defined on a subset $\Omega$ of
  the unit disk $\boldD$ with values in $\eL(\eG,\eK)$ and
  $\eL(\eF,\eK)$, where $\eF,\eG,\eK$ are Hilbert spaces.  Assume that
  both
  $$
  K_{1}(w,z) = \dfrac{A(z)A(w)^*-B(z)B(w)^*}{1-z\bar w}
  $$
  and
  $$
  K_{2}(w,z) = \dfrac{A(z)A(w)^*-z\bar w B(z)B(w)^*}{1-z\bar w}
  $$
  have $\kappa$ negative squares on $\Omega \times \Omega$.  Then
  there is a function $S(z)$ in $\SkappaprimeFG$, $\kappa' \le
  \kappa$, such that $B(z) = A(z)S(z)$ for all but at most $\kappa$
  points $z$ of $\Omega$.  In this case, $\kappa' = \kappa$ if and
  only if the elements $h$ of $\eH(S)$ such that $A(z)h(z) \equiv 0$
  on $\Omega$ form a Hilbert subspace of $\eH(S)$.
\end{theorem}

The proof uses a different colligation from that of
Theorem~\ref{rectangular}.  It is adapted from the work of
V.~E.~Katsnelson, A.  Kheifets, and P.~M. Yuditski\u\i; see Kheifets
\cite{Kheifets98} for an account and references to earlier works.  The
idea is used by Ball and Trent \cite{BT}, who extend it to a several
variable setting and apply it in a form for reproducing kernel
functions that is close to our situation.

Theorem~\ref{th3} is a non-holomorphic analog of \cite[Theorem
11]{ADRSlivsic}: there the coefficient spaces are indefinite, but we
have the stronger hypothesis that $\Omega$ is a neighborhood of the
origin and $A(z)$ and $B(z)$ are holomorphic.  Now the functions
$A(z)$ and $B(z)$ are not assumed to be holomorphic, but in
compensation $\eF$ and $\eG$ are required to be Hilbert spaces (for
simplicity we have taken $\eK$ to be a Hilbert space also, but this
plays no role in the argument).  The proof of Theorem~\ref{th3} runs
along the same lines.

\begin{proof}
  Write $\eH(K_{1})$ and $\eH(K_{2})$ for the Pontryagin spaces with
  reproducing kernels $K_{1}(w,z)$ and $K_{2}(w,z)$.  Define a
  relation
  \begin{equation}\notag
    \boldR = \text{span} \left\{
\left(
\left(
\begin{matrix}
  K_{1}(w,\cdot)k \\ B(w)^*k
\end{matrix}
\right),
\left(
\begin{matrix}
  \bar w K_{1}(w,\cdot)k \\ A(w)^*k
\end{matrix}
\right)
\right) : w\in\Omega,\; k\in\eK
\right\}
\subseteq \left(
\begin{matrix}
  \eH(K_{1}) \\ \eF
\end{matrix}
\right)
\times
\left(
\begin{matrix}
  \eH(K_{1}) \\ \eG
\end{matrix}
\right).
  \end{equation}
  It is easy to see that $\boldR$ is isometric.  We show that the
  domain $\eM$ of $\boldR$ contains a maximal uniformly negative
  subspace of $\eH(K_{1})\oplus \eF$.  To this end, consider a Gram
  matrix of the form
  $$
  M= \left( {\ip {\left(
\begin{matrix}
  K_{1}(w_j,\cdot)k_j \\ B(w_j)^*k_j
\end{matrix}
\right)} {\left(
\begin{matrix}
  K_{1}(w_i,\cdot)k_i \\ B(w_i)^*k_i
\end{matrix}
\right)} }_{\eH(K_{1})\oplus \eF} \right)_{i,j=1}^n ,
$$
where $w_1,\dots,w_n$ are any points in $\Omega$ and
$k_1,\dots,k_n$ are arbitrary vectors in $\eK$.  Thus
\begin{equation}\notag
  M = \Big(
{\ip{\,\left[ K_{1}(w_j,w_i)+ B(w_i)B(w_j)^*\right]k_j}{k_i}}_{\eK}
\Big)_{i,j=1}^n 
 = \Big(
{\ip{K_{2}(w_j,w_i)k_j}{k_i}}_{\eK}
\Big)_{i,j=1}^n .
\end{equation}
Since we assume that $K_{2}(w,z)$ has $\kappa $ negative squares, $M$
has at most $\kappa $ negative eigenvalues no matter how
$w_1,\dots,w_n$ and $k_1,\dots,k_n$ are chosen, and some such Gram
matrix has exactly $\kappa $ negative eigenvalues.  By \cite[Lemma
$1.1.1'$]{ADRSbook}, $\eM$ contains a $\kappa$-dimensional subspace
which is the antispace of a Hilbert space in the inner product of
$\eH(K_{1})\oplus \eF$.  Since $\sqminus (\eH(K_{1})\oplus \eF) =
\kappa $, this verifies the assertion.  It follows that the closure of
$\eM$ in $\eH(K_{1})\oplus \eF$ is a regular subspace whose orthogonal
complement {\smash{$\eM^\perp$}} is a Hilbert space.
  
By \cite[Theorem 1.4.2]{ADRSbook}, the closure of the range of
$\boldR$ is likewise a regular subspace $\eN$ of $\eH(K_{1})\oplus
\eG$, and we can construct a partial isometry
$$
V = \left(\begin{matrix} T&F\\G&H \end{matrix}\right) \colon
\begin{pmatrix}
  \eH(K_{1}) \\ \eF
\end{pmatrix}
\to
\begin{pmatrix}
  \eH(K_{1}) \\ \eG
\end{pmatrix}
$$
with initial space $\eM$ and final space $\eN$ such that
$$
V^*= \left(\begin{matrix} T^*&G^*\\F^*&H^* \end{matrix}\right) :
\left( \begin{matrix} \bar w K_{1}(w,\cdot)k \\ \noalign{\vskip4pt}
    A(w)^*k
  \end{matrix} \right) \rightarrow \left( \begin{matrix}
    K_{1}(w,\cdot)k \\ \noalign{\vskip4pt} B(w)^*k \end{matrix}
\right)
$$
for all $k\in\eK$ and all $w\in\Omega$.  Thus for $w\in\Omega$,
\begin{equation}\label{E:id1}
T^*\left\{ \bar w K_{1}(w,\cdot)k \right\} + 
G^*\left\{ A(w)^*k \right\} = K_{1}(w,\cdot)k,
\end{equation}
and
\begin{equation}\label{E:id2}
F^*\left\{ \bar w K_{1}(w,\cdot)k \right\} + 
H^*\left\{ A(w)^*k \right\} =B(w)^*k.
\end{equation}
Hence
\begin{equation}\label{E:id3}
(1-\bar w T^*)\left\{ K_{1}(w,\cdot) A(w)^*k\right\} 
= G^*\left\{ A(w)^*k \right\}.
\end{equation}
Since $\ker V$ is a Hilbert space, $V$ is a contraction.  As in the
proof of Theorem \ref{rectangular}, because we assume that $\eF$ and
$\eG$ are Hilbert spaces, $T$ is a contraction, and the part of the
spectrum of $T$ that lies in $|\lambda|>1$ consists of at most
$\kappa$ normal eigenvalues.

Let $\Omega' = \Omega \setminus \{ \lambda_1,\dots,\lambda_q \}$,
where $\lambda_1,\dots,\lambda_q$ are the points $\lambda$ of the unit
disk at which $1-\lambda T$ is not invertible ($q\le\kappa$).  For all
$w\in\Omega'$ and all $k\in\eK$,
$$
K_{1}(w,\cdot)k = (1-\bar w T^*)^{-1}G^* \{A(w)^*k\}
$$
by \eqref{E:id3}.  Define
$$
S(z) = H + zG(1-zT)^{-1}F, \qquad z\in \boldD \setminus \{
\lambda_1,\dots,\lambda_q \}.
$$
Then $B(w) = A(w)S(w)$, $w\in \Omega'$, by \eqref{E:id1} and
\eqref{E:id2}.  The proof that $\hat S\in\Skappaprime$ for some
$\kappa'\le\kappa$ is the same as in the proof of
Theorem~\ref{rectangular}.  The last statement follows from
\cite[Theorem 1.5.7]{ADRSbook}.
\end{proof}

In the next theorem, we allow $\eF,\eG,\eK$ to be indefinite, but the
functions $A(z)$ and $B(z)$ are required to be holomorphic.  This
yields a new result of Leech type factorization theorems as a
companion to those of \cite{ADRSlivsic}.

\begin{theorem}\label{rectangular2}
  Let $\eF$, $\eG$, $\eK$ be Kre{\u\i}n spaces with $\sqminus \eF =
  \sqminus \eG < \infty$.  Let $\Omega$ be a subregion of the unit
  disk containing the origin.  Let $A(z)$ and $B(z)$ be holomorphic
  functions on $\Omega$ with values in $\eL(\eG,\eK)$ and
  $\eL(\eF,\eK)$.  Assume that the kernel
  \begin{equation}
    \label{ABkernel3}
    K(w,z) = \frac{A(z)A(w)^* - B(z)B(w)^*}{1 - \bar w z}
  \end{equation}
  has $\kappa$ negative squares on $\Omega\times\Omega$, and let
  $\eH_K$ be the associated reproducing kernel Pontryagin space.  Let
  $\eM$ be the subspace of $\eH_K\oplus \eG$ consisting of all
  elements $k(z)\oplus g$ such that
  $$
  A(0)g=0 \quad \text{and}\quad zk(z) + [A(z)-A(0)]g \equiv 0 \quad
  \text{on}\quad \Omega.
  $$
  Let $\eN$ be the subspace of $\eH_K\oplus \eF$ consisting of all
  elements $h(z)\oplus f$ such that
  $$
  h(z) + B(z)f \equiv 0 \quad \text{on}\quad \Omega.
  $$
  Assume that $\eM$ and $\eN$ are Hilbert spaces in the inner
  products of the larger spaces.  Then there is a function $S(z) \in
  \SkappaprimeFG$ for some $\kappa'\le\kappa$ which is holomorphic at
  the origin and such that $B(z) = A(z)S(z)$ for all but at most
  $\kappa$ points $z$ of $\Omega$.  In this case, $\kappa' = \kappa$
  if and only if the elements $h$ of $\eH(S)$ such that $A(z)h(z)
  \equiv 0$ on $\Omega$ form a Hilbert subspace of $\eH(S)$.
\end{theorem}

\begin{proof}
  We repeat the constructions in the proof of
  Theorem~\ref{rectangular}.  The partial isometry $V$ is again a
  contraction in the present situation.  In general, the operator $T$
  is not a contraction, but it is a bounded operator and so
  $(1-wT)^{-1}$ is defined for $|w|$ sufficiently small.  The argument
  goes through if we restrict attention to a suitable neighborhood of
  the origin.  At the end, the identity $B(z) = A(z)S(z)$ extends to
  all but at most $\kappa$ points of $\Omega$ by analytic
  continuation.
\end{proof}

\section{Coefficient and moment problems}\label{S:coefficients}

Let $z_1,\dots,z_n$ be points in the unit disk, and let
$w_1,\dots,w_n$ be any complex numbers.  If we specialize
Section~\ref{S:interpolation} to the scalar case and set $\Omega=\{
z_1,\dots,z_n\}$, $A(z_j)=1$, and $B(z_j)=w_j$ for all $j=1,\dots,n$,
then the interpolation problem in Section~\ref{S:interpolation}
reduces to the Nevanlinna-Pick problem.  The indefinite form of
interpolation was introduced by Takagi \cite{Takagi}, and it has been
studied by Adamjan, Arov, and Kre{\u\i}n \cite{MR45:7505}, Kre{\u\i}n
and Langer \cite{MR57:1173}, and others.  A rather complete picture of
the solution of the indefinite Nevanlinna-Pick problem emerged from
this work.  A remaining issue concerning the degenerate case was
recently settled.  Namely, one can ask, for which nonnegative integers
$\kappa$ can the Nevanlinna-Pick problem be solved in ${\bf
  S}_{\kappa}$ for given data $z_1,\dots,z_n$ and $w_1,\dots,w_n$?  A
more precise question can be posed.  Define ${\bf S}_{\nu, \pi}$ as
the class of all meromorphic functions $S(z)$ on the unit disk for
which the kernel $K_S(w,z)$ has $\nu$ negative squares and $\pi$
positive squares (thus ${\bf S}_{\nu, \pi}$ is a subclass of ${\bf
  S}_{\nu}$).  For which nonnegative integers $\nu$ and $\pi$ can the
Nevanlinna-Pick problem be solved in ${\bf S}_{\nu, \pi}$ for given
data $z_1,\dots,z_n$ and $w_1,\dots,w_n$?  These questions were
answered by Woracek \cite{MR96k:47021, MR98e:30041} (with the disk
replaced by the upper half-plane), yielding a complete solution of the
Nevanlinna-Pick problem in the scalar case.

We consider analogous questions for the indefinite
Carath\'eodory-Fej\'er problem and obtain a complete solution in the
scalar case.  The solution depends on results of Iokhvidov
\cite{Iokhvidov} on a related trigonometric moment problem.  In the
positive definite case this connection is well-known.  We refer to
\cite{Iokhvidov} for references to the original papers (some jointly
with M.~G. Kre{\u\i}n) pertaining to this problem.  A key step
involves another application of the characteristic function of a
partially isometric operator colligation, which was the principal tool
in Section~\ref{S:interpolation}.

\begin{problem}[Carath\'eodory-Fej\'er problem]\label{P:CF} 
  Let $a_0,a_1,\dots,a_{n-1}$ be $n$ complex numbers. For which
  nonnegative integers $\kappa$ is there a function $S(z)$ in
  $\Skappa$ which is holomorphic at the origin and such that $S(z) =
  a_0 + a_1 z + \cdots + a_{n-1}z^{n-1} + \oh(z^n)$ in a neighborhood
  of the origin?  For which $\nu$ and $\pi$ do there exist solutions
  in ${\bf S}_{\nu, \pi}$?

\end{problem}

Necessary conditions on coefficients are obtained from the series
expansions of standard kernel functions.  Suppose that $S(z)$ is a
holomorphic (scalar-valued) function defined in a neighborhood of the
origin.  Let $S(z) = a_0 +a_1 z + a_2 z^2 + \cdots$ be its Taylor
series expansion, and write
\begin{align}
  T_r &= \left(
\begin{matrix}
  a_0 & 0 & 0 & \cdots & 0 \\
  a_1 & a_0 & 0 & \cdots  & 0\\
  && \cdots && \\
  a_{r-1} & a_{r-2} & a_{r-3} & \cdots  & a_0\\
\end{matrix}
\right), \qquad \tilde T_r = \left(
\begin{matrix}
  \bar a_0 & 0 & 0 & \cdots & 0 \\
  \bar a_1 & \bar a_0 & 0 & \cdots  & 0\\
  && \cdots && \\
  \bar a_{r-1} & \bar a_{r-2} & \bar a_{r-3} & \cdots  & \bar a_0\\
\end{matrix}
\right),
                              \label{E:toeplitz1}\\
                              &\hskip2.5cm Q_r = \left(
\begin{matrix}
  a_1 & a_2 & a_3 & \cdots & a_{r} \\
  a_2 & a_3 & a_4 & \cdots  & a_{r+1}\\
  && \cdots && \\
  a_{r} & a_{r+1} & a_{r+2} & \cdots & a_{2r-1} \\
\end{matrix}  
\right), \label{E:hankel1}
\end{align}
$r=1,2,\dots$.  Set $\tilde S(z) = \bar S(\bar z)$.  Straightforward
calculations yield the expansions {\allowdisplaybreaks
\begin{align}
  K_{S}(w,z) &= \frac{1-S(z)\overline{S(w)}}{1-z\bar w} =
  \sum_{p,q=0}^\infty\, C_{pq}z^p\bar w^q,  \notag \\
  K_{\tilde S}(w,z) &= \frac{1-\tilde S(z)\overline{\tilde
      S(w)}}{1-z\bar w} = \sum_{p,q=0}^\infty\,
  \tilde C_{pq}z^p\bar w^q, \notag \\
  D_S(w,z) &= \left(
\begin{matrix}
  K_S(w,z) & \displaystyle{\frac{S(z)-S(\bar w)}{z-\bar w}} \\
  \displaystyle{\frac{\tilde S(z)- \tilde S(\bar w)}{z-\bar w}} & K_{\tilde S}(w,z) \\
\end{matrix}
\right) = \sum_{p,q=0}^\infty\, D_{pq}z^p\bar w^q, \notag
\end{align}
} $\!\!$where {\allowdisplaybreaks\begin{align} \big[ C_{pq}
  \big]_{p,q=0}^{n-1} &= I_n - T_nT_n^*, \hskip1.75cm \big[ \tilde
  C_{pq} \big]_{p,q=0}^{n-1} = I_n - \tilde T_n \tilde T_n^*,
  \notag\\
  &\hskip-1.5cm \big[ D_{pq} \big]_{p,q=0}^{r-1} =
  \left(\begin{matrix}
      I_r - T_rT_r^*  & Q_r \\
      \noalign{\vskip4pt} Q_r^* & I_r - \tilde T_r \tilde T_r^*
\end{matrix}\right),  \qquad 1 \le r \le n/2. \notag
\end{align}
} $\!$Thus the coefficients $a_0,a_1, \dots$ of $S(z)$ give rise to
three families of matrices:
\begin{equation}\label{E:toeplitz2}
I_n - T_nT_n^*, \qquad I_n - \tilde T_n \tilde T_n^*,\qquad
\left(\begin{matrix}
I_r - T_rT_r^*  & Q_r \\
\noalign{\vskip4pt}
Q_r^*     & I_r - \tilde T_r \tilde T_r^*
\end{matrix}\right),
\qquad 1 \le r \le n/2,
\end{equation}
$n=0,1,\dots$.  For fixed $n$, the matrices \eqref{E:toeplitz2} depend
only on $a_0,\dots,a_{n-1}$.

If $S(z)$ belongs to $\Skappa$, then the three kernels each have
$\kappa$ negative squares \cite[Theorem 2.5.2]{ADRSbook}.  It follows
that the number of negative eigenvalues of each of the matrices in
\eqref{E:toeplitz2} is a nondecreasing function of the order of the
matrix, and this number is ultimately equal to~$\kappa$ in each case
(see the result in Section~\ref{S:kernels}).

If $S(z)$ belongs to ${\bf S}_{\nu, \pi}$, similar remarks apply not
only to the number of negative squares but also to the number of
positive squares.  For simplicity, suppose that $S(0) \neq 0$, and
note the identities
\begin{align}
  K_S(w,z) &= - S(z) K_{1/S}(w,z) \overline{S(w)}, \notag \\
  K_{\tilde S}(w,z) &= - \tilde S(z) K_{1/{\tilde S}}(w,z)
  \overline{\tilde S(w)}, \notag \\
  D_S(w,z) &= - \begin{pmatrix} S(z) & 0 \\ 0 & \tilde S(z)
  \end{pmatrix}
  D_{1/S}(w,z)
  \begin{pmatrix}
    \overline{S(w)} & 0 \\ 0 & \overline{\tilde S(w)}
  \end{pmatrix} . \notag
\end{align}
The numbers of positive squares of $K_{S}(w,z)$, $K_{\tilde S}(w,z)$,
and $D_{S}(w,z)$ thus coincide, with the numbers of negative squares
of $K_{1/S}(w,z)$, \smash{$K_{1/\tilde S}(w,z)$,} and $D_{1/S}(w,z)$,
respectively.  Thus if one of the three kernels has $\pi$ positive
squares, then all do.  In this case, applying the previous assertions
concerning negative squares, we see that the number of positive
eigenvalues of each of the matrices in \eqref{E:toeplitz2} is a
nondecreasing function of the order of the matrix, and this number is
ultimately equal to~$\pi$ in each case.

This raises questions concerning the general behavior of the numbers
of negative and positive eigenvalues for the matrices
\eqref{E:toeplitz2} whenever \eqref{E:toeplitz1} and \eqref{E:hankel1}
are defined for any complex numbers $a_0,a_1,\dots$, whether these
numbers are the Taylor coefficients of a holomorphic function or not.
We show that the behavior is indeed always similar to the special
cases noted above: the numbers of negative (positive) eigenvalues for
the three types are nondecreasing functions of the order, and if one
eventually has some constant value, then all have the same constant
value eventually.  These questions are purely algebraic.  There is a
separate convergence question, namely, under what conditions are the
given numbers $a_0,a_1,\dots$ the Taylor coefficients of a holomorphic
function $S(z)$ in $\Skappa$ or ${\bf S}_{\nu, \pi}$?  Finally, if we
only define \eqref{E:toeplitz1}, \eqref{E:hankel1}, and
\eqref{E:toeplitz2} as far as we can go with a finite sequence
$a_0,\dots,a_{n-1}$, what are the possible extensions to an infinite
sequence $a_0,a_1,\dots$?

To answer such questions, we relate given complex numbers
$a_0,\dots,a_{n-1}$ to a trigonometric moment problem.  Define
$c_0=1,c_1, \dots,c_n$ by
\begin{equation}\label{E:connection}
\left\{
\begin{array}{cl}
c_0 \!\!\! &= 1,                      \\
c_1 \!\!\! &= c_0 a_0,                \\
c_2 \!\!\! &= c_0 a_1 + c_1 a_0,       \\
 \!\!\! & \dots                     \\
c_n \!\!\! &= c_0 a_{n-1} + c_1 a_{n-2} + \cdots + c_{n-1} a_0,
\end{array}  
\right.
\end{equation}
This correspondence is one-to-one and has the property that if
$a_0,\dots,a_{n-1}$ corresponds to $c_0=1,c_1, \dots,c_n $ then for
each $1 \le k \le n$, $a_0,\dots,a_{k-1}$ corresponds to
$c_0=1,c_1,\dots,c_k $ also via (\ref{E:connection}) with $n$ replaced
by $k$.  We consider the associated matrix
\begin{equation}\label{E:tmp1}
    M_n = \begin{pmatrix}
   c_0 & \bar c_1 & \bar c_2 & \dots & \bar c_n     \\[3pt]
   c_1 & c_0      & \bar c_1 & \dots & \bar c_{n-1} \\[3pt]
   c_2 & c_1      & c_0      & \dots & \bar c_{n-2} \\[3pt]
       &          & \dots                      \\[3pt]
   c_n & c_{n-1}  & c_{n-2}  & \dots & c_0    
    \end{pmatrix}
\end{equation}
In the sequel $J_n$ stands for the selfadjoint and unitary $n \times
n$ matrix
$$
J_n = \left(\begin{matrix}
    0 & 0 &  \dots & 0 & 1 \\
    0 & 0 &  \dots & 1 & 0 \\
    && \dots && \\
    1 & 0 & \dots & 0 & 0 \\
\end{matrix}
\right).
$$
Also define
$$
B_r = \left(\begin{matrix}
    c_0 & 0 &  \dots & 0  \\
    c_1 & c_0 &  \dots  & 0 \\
    && \dots & \\
    c_r & c_{r-1} & \dots & c_0 \\
\end{matrix}
\right), \quad C_r = \left(\begin{matrix}
    I_r & 0 \\
    0 & B_r \\
\end{matrix}
\right) \left(\begin{matrix}
    0 & B_r^*J_{r+1} \\
    I_r & 0 \\
\end{matrix}
\right).
$$

\begin{theorem} \label{matrixident} 
  Let $a_0,a_1,\dots,a_{n-1}$ be complex numbers and define
  $c_0=1,c_1,\dots,c_n$ by \eqref{E:connection}.  The following
  equalities hold:
\begin{equation}\label{E:stronger}
M_r = B_r 
\begin{pmatrix}
  1  & 0  \\
  0 & I_r - T_rT_r^*
\end{pmatrix}
B_r^*
= B_r^* 
  \begin{pmatrix}
    I_r-T_r^*T_r &0 \\
    0 & 1
  \end{pmatrix} B_r, \quad  1 \le r \le n,
\end{equation}
\begin{equation}\label{E:barmr}
\overline{M}_r =\tilde B_r^*
\begin{pmatrix} I_r-\tilde T_r^* \tilde T_r & 0 \\
  0 & 1 \end{pmatrix} \tilde B_r, \quad 1\le r \le n,
\end{equation}
and
\begin{equation}\label{E:bigmatrix}
M_{2r} = C_r  
\begin{pmatrix}
  I_r-T_rT_r^* &0 &Q_r \\
  0&1&0\\
  Q_r^* &0&I_r-\tilde{T}_r \tilde{T}_r^*
\end{pmatrix} C_r^*, \quad 1 \le r \le n/2.
\end{equation}
\end{theorem}

The bar in $\overline{M}_r$ in (\ref{E:barmr}) indicates that all
entries in the matrix $M_r$ have been replaced by their complex
conjugates.

\begin{proof}
  The first equality in (\ref{E:stronger}) can be shown by induction.
  The second equality follows from the first.  To see this, use the
  identities $J_{r+1}M_rJ_{r+1}=\overline{M}_r$,
  $J_{r+1}B_rJ_{r+1}=\tilde B _r^*$, $J_rT_rJ_r = \tilde T_r^*$, and
  $$\begin{pmatrix} 0 & J_r \\ 1 & 0 \end{pmatrix}
\begin{pmatrix} 1 & 0 \\ 0 & I_r-T_rT_r^* \end{pmatrix}
\begin{pmatrix} 0 & 1 \\ J_r & 0 \end{pmatrix}
=
\begin{pmatrix} I_r-\tilde T_r^* \tilde T_r & 0 \\
  0 & 1 \end{pmatrix},
$$
to obtain
$$
M_r = J_{r+1}\overline{M}_rJ_{r+1} = J_{r+1}\overline{B}_rJ_{r+1}
\begin{pmatrix}
  I_r - T_r^*T_r & 0 \\
  0 & 1
\end{pmatrix}
J_{r+1}\overline{B}_r^*J_{r+1} = B_r^* \begin{pmatrix}
  I_r - T_r^*T_r & 0 \\
  0 & 1
\end{pmatrix}
B_r,
$$
which is the second equality in \eqref{E:stronger}.  We get
\eqref{E:barmr} on replacing the entries of the matrices by their
complex conjugates.

We prove (\ref{E:bigmatrix}). Assume $1 \le r \le n/2$.  Then
  \begin{equation}\label{E:twotimesr}
    M_{2r} = 
    \begin{pmatrix}
      M_{r-1} & S_r^* \\ S_r & M_r
    \end{pmatrix},
\qquad
 S_r = 
   \begin{pmatrix}
     c_{r} & c_{r-1} & \dots &c_{2} &c_{1} \\
     c_{r+1} & c_{r} & \dots &c_{3} &c_{2} \\
     &&& \dots &  \\
     c_{2r} & c_{2r-1} & \dots &c_{r+2} &c_{r+1} \\
   \end{pmatrix}.
  \end{equation}
  In \eqref{E:twotimesr} we use the first equality in
  \eqref{E:stronger} to obtain
\begin{multline}\notag
  M_{2r} =
  \begin{pmatrix}
    M_{r-1} & S_r^* \\
    S_r & B_r
      \begin{pmatrix} 1 & 0 \\ 0 & I_r-T_rT_r^*\end{pmatrix} B_r^*
  \end{pmatrix}  \\
  =
\begin{pmatrix}
  I_r&0\\0&B_r  \\
\end{pmatrix}
\begin{pmatrix}
  M_{r-1} & S_r^*B_r^{*-1}\\
  B_r^{-1}S_r & \begin{pmatrix} 1 & 0 \\ 0 & I_r-T_rT_r^*\end{pmatrix}
\end{pmatrix}
\begin{pmatrix}
  I_r & 0 \\ 0 & B_r^*
\end{pmatrix}.
\end{multline}
Due to the lower triangular form of $B_r$, we get
\begin{equation}\notag
  B_r^{-1} S_r = 
  \begin{pmatrix}
    1 &0&\dots &0&0 \\
    * &1&\dots &0&0 \\
    && \dots && \\
    *&*&\dots &*&1
  \end{pmatrix}
  \begin{pmatrix}
    c_{r}&c_{r-1}&\dots &c_{1} \\
    &&& \\
    && \hskip-1cm * & \\
    &&& \\
  \end{pmatrix}
= 
\begin{pmatrix}
  c_{r}&c_{r-1}&\dots &c_{1} \\
  &&& \\
  && \hskip-1cm Z_r & \\
  &&& \\
  \end{pmatrix}.
\end{equation}
With this definition of $Z_r$ and (\ref{E:stronger}), we obtain
{\allowdisplaybreaks
\begin{align}
  M_{2r} &= \begin{pmatrix} I_r&0\\[12pt]0&B_r
  \end{pmatrix} 
\begin{pmatrix}
  M_{r-1} &
    \begin{pmatrix}
      \bar c_r && \; & \; & \; &\\
      \vdots && \; & Z_r^* & \;& \\
      \bar c_1 && \; & \; & \; & \;& \\
    \end{pmatrix}
    \\[25pt]
    \begin{pmatrix}
      c_{r}&c_{r-1}&\dots &c_{1} \\
      &&& \\
      && \hskip-1cm Z_r \vphantom{\vdots} & \\
      &&& \\
  \end{pmatrix} 
  &
\begin{pmatrix}
  1 & 0 & \dots & 0 \\
  0 & \; & \; & \; \\
  \vdots & \; & I_r-T_rT_r^* & \; \\
  0& \; & \; & \; \\
\end{pmatrix}
  \end{pmatrix}
  \begin{pmatrix}
    I_r&0\\[12pt]0&B_r^*
  \end{pmatrix}
  \notag \\[9pt]
  &= \begin{pmatrix} I_r&0\\[12pt]0&B_r
  \end{pmatrix}
  \begin{pmatrix}
    M_r & \begin{pmatrix}
      && \\
      & Z_r^* &  \\
      && \\
      0 &\dots &0
\end{pmatrix} \\ 
\begin{pmatrix}
  &&&& 0 \\
  && Z_r &&\vdots \\
  &&&& 0 \end{pmatrix} &
\begin{matrix}
  && \\
  &I_r-T_rT_r^*& \\
  && \\
\end{matrix}
\end{pmatrix}
\begin{pmatrix}
  I_r&0\\[12pt]0&B_r^*
  \end{pmatrix}
  \notag \\[9pt]
  &= \begin{pmatrix} I_r&0\\[12pt]0&B_r
  \end{pmatrix}
  \begin{pmatrix}
    B_r^* \begin{pmatrix} I_r-T_r^*T_r & 0 \\ 0&1
\end{pmatrix}B_r & \begin{pmatrix}
  && \\
  & Z_r^* &  \\
  && \\
  0 &\dots &0
\end{pmatrix} \\ 
\begin{pmatrix}
  &&&& 0 \\
  && Z_r &&\vdots \\
  &&&& 0 \end{pmatrix} &
\begin{matrix}
  && \\
  &I_r-T_rT_r^*& \\
  && \\
\end{matrix}
  \end{pmatrix}
\begin{pmatrix}
  I_r&0\\[12pt]0&B_r^*
  \end{pmatrix} \notag \\[9pt]
  &=
\begin{pmatrix}
  I_r&0\\[12pt]0&B_r
  \end{pmatrix}
\begin{pmatrix}
  B_r^*&0\\[12pt]0&I_r
  \end{pmatrix}
  \begin{pmatrix}
    \begin{pmatrix}
      I_r-T_r^*T_r &0 \\0 & 1
    \end{pmatrix}
    & B_r^{*-1}
\begin{pmatrix}
  &&\\
  & Z_r^*& \\
  &&& \\
  0&\dots& 0
\end{pmatrix}
\\[25pt]
\begin{pmatrix}
  &&&0\\
  &Z_r &&\vdots \\
  &&&0
\end{pmatrix}B_r^{-1}
& I_r-T_rT_r^*
  \end{pmatrix}
  \cdot \notag \\[17pt]
  &\hskip6cm \cdot
\begin{pmatrix}
  B_r&0\\[12pt]0&I_r
  \end{pmatrix}
\begin{pmatrix}
  I_r & 0 \\[12pt] 0&B_r^*
\end{pmatrix} \notag
\end{align}
} 
Here the matrix
$$
C_r' = \begin{pmatrix} I_r&0\\[12pt]0&B_r
  \end{pmatrix}
\begin{pmatrix}
  B_r^*&0\\[12pt]0&I_r
  \end{pmatrix}
  $$
  is invertible.  Note also that
  $$
\begin{pmatrix}
  &&&0\\
  &Z_r &&\vdots \\
  &&&0
\end{pmatrix}B_r^{-1} = \begin{pmatrix}
  &&&0\\
  &Z_r &&\vdots \\
  &&&0
\end{pmatrix} 
\begin{pmatrix}
  1&0 &\dots &0 \\
  *&1 &\dots &0 \\
  & \dots && \\
  *&*&\dots&1
\end{pmatrix}
=
\begin{pmatrix}
  &&&0\\
  &Y_r &&\vdots \\
  &&&0
\end{pmatrix},
$$
so that with $Y_r$ defined in this way, we have
\begin{equation}\label{E:almostthere}
M_{2r} = C_r'  
\begin{pmatrix}
  I_r-T_r^*T_r &0 &Y_r^* \\
  0&1&0\\
  Y_r &0&I_r- T_r T_r^*
\end{pmatrix} {C_r '}^*, \quad
1 \le r \le n/2.
\end{equation}
We now identify $Y_r$ as
\begin{equation}\label{E:Ysubr}
  Y_r = 
  \begin{pmatrix}
    a_{r} & a_{r-1} &\dots & a_{1} \\
    a_{r+1} & a_{r} &\dots & a_{2} \\
    & \hskip.5cm\dots &&  \\
    a_{2r-1} & a_{2r-2} &\dots & a_{r} \\
  \end{pmatrix} .
\end{equation}
>From the definition of $B_r$ we find that
\begin{equation}\notag
  B_r^{-1} = 
  \begin{pmatrix}
    1 & 0 & \dots & 0 \\
    -a_0 & 1 & \dots & 0 \\
    & \hskip.5cm\dots &&  \\
    -a_{r-1} & -a_{r-2} & \dots & 1 \\
  \end{pmatrix}.
\end{equation}
It follows that
\begin{equation}\notag
  B_r^{-1}S_r = 
  \begin{pmatrix}
    \begin{pmatrix} c_r & c_{r-1} & \dots & c_1 \end{pmatrix}
    \\[5pt]
  \begin{pmatrix}
    a_{r} & a_{r-1} &\dots & a_{1} \\
    a_{r+1} & a_{r} &\dots & a_{2} \\
    & \hskip.5cm\dots &&  \\
    a_{2r-1} & a_{2r-2} &\dots & a_{r} \\
  \end{pmatrix}B_{r-1}
  \end{pmatrix}
\quad  {\rm and} \quad Z_r = \begin{pmatrix}
    a_{r} & \!a_{r-1} &\dots & a_{1} \\
    a_{r+1} \!& a_{r} &\dots & a_{2} \\
& \hskip.5cm\dots \!&&  \\
    a_{2r-1} \!& a_{2r-2} &\dots & a_{r} \\
  \end{pmatrix}B_{r-1}.
\end{equation}
Finally, we obtain
\begin{align}\notag
\begin{pmatrix} Y_r & 0 \end{pmatrix}
=
  \begin{pmatrix} Z_r & 0 \end{pmatrix} B_r^{-1}
  &=\begin{pmatrix}
  \begin{pmatrix}
    a_{r} & a_{r-1} &\dots & a_{1} \\
    a_{r+1} & a_{r} &\dots & a_{2} \\
    & \hskip.5cm\dots &&  \\
    a_{2r-1} & a_{2r-2} &\dots & a_{r} \\
  \end{pmatrix}B_{r-1} & 0
\end{pmatrix}
\begin{pmatrix}
  B_{r-1}^{-1} & 0 \\[12pt] * &1
\end{pmatrix}
\notag \\
&\hskip1cm =
\begin{pmatrix}
  \begin{pmatrix}
    a_{r} & a_{r-1} &\dots & a_{1} \\
    a_{r+1} & a_{r} &\dots & a_{2} \\
    & \hskip.5cm\dots &&  \\
    a_{2r-1} & a_{2r-2} &\dots & a_{r} \\
  \end{pmatrix} &0
\end{pmatrix},
\notag
\end{align}
proving \eqref{E:Ysubr}. Evidently, $Q_r=Y_rJ_r$ and
$$C_r=C'_r\begin{pmatrix} 0 & J_{r+1} \\ I_r & 0 \end{pmatrix}
=C'_r\begin{pmatrix} 0 & 0 & J_r \\ 0 & 1 & 0 \\ I_r & 0 & 0
\end{pmatrix} .
$$
Substituting this in (\ref{E:almostthere}) we obtain
(\ref{E:bigmatrix}).
\end{proof}

A number of consequences follow. For any Hermitian matrix $A$ we write
$\pi(A)$ and $\nu(A)$ for the numbers of positive and negative
eigenvalues of $A$ counting multiplicity.

\begin{corollary}
  Let $a_0,a_1,\dots,a_{n-1}$ be complex numbers and define
  $c_0=1,c_1,\dots,c_n$ by \eqref{E:connection}.
  
  $(1)$ Each of the four quantities
\begin{gather}
  \nu (I_r - T_r T_r^*), \quad \pi (I_r - T_r T_r^*),
  \qquad 1 \le r \le n, \notag \\[6pt]
  \nu \begin{pmatrix} I_r-T_r T_r^* &Q_r \\ Q_r^*&I_r-\tilde T_r\tilde
    T_r^*
\end{pmatrix} , \quad 
\pi \begin{pmatrix} I_r-T_r T_r^* &Q_r \\ Q_r^*&I_r-\tilde T_r\tilde
  T_r^*
\end{pmatrix},  
\qquad 1 \le r \le n/2, \notag
\end{gather}
is a nondecreasing function of $r$.

$(2)$ For $0 \le r \le n$,
\begin{equation}\notag
  \nu (I_r - T_r T_r^*)
   = \nu (I_r - T_r^*T_r)
   = \nu (I_r - \tilde T_r \tilde T_r^*)
   = \nu (I_r - \tilde T_r^*\tilde T_r)
\end{equation}
and
\begin{equation}\notag
  \pi (I_r - T_r T_r^*)
   = \pi (I_r - T_r^*T_r)
   = \pi (I_r - \tilde T_r \tilde T_r^*)
   = \pi (I_r - \tilde T_r^*\tilde T_r).
\end{equation}
$(3)$ If $\nu (I_n - T_nT_n^*) = \kappa$, then all of the matrices in
$\eqref{E:toeplitz2}$ have at most $\kappa$ negative eigenvalues.
\end{corollary}

The condition $\nu (I_n - T_nT_n^*) = \kappa$ is necessary that
$a_0,a_1,\dots,a_{n-1}$ are the first $n$ Taylor coefficients of a
function in $\Skappa$.  The point of statement (3) in the preceding
corollary is that no stronger necessary condition can be obtained from
the other matrices in \eqref{E:toeplitz2}.

\begin{proof}
  (1) By the first equality in \eqref{E:stronger},
  \begin{align}
    \nu (I_r - T_r T_r^*) &= \nu (M_r), \notag \\
    \pi (I_r - T_r T_r^*) &= \pi (M_r) - 1. \notag
  \end{align}
  By \eqref{E:bigmatrix},
  \begin{align}
    \nu \begin{pmatrix} I_r-T_r T_r^* &Q_r \\ Q_r^*&I_r-\tilde
      T_r\tilde T_r^*
\end{pmatrix}&= \nu (M_{2r})  , \notag \\
\pi \begin{pmatrix} I_r-T_r T_r^* &Q_r \\ Q_r^*&I_r-\tilde T_r\tilde
  T_r^*
\end{pmatrix} &= \pi (M_{2r}) - 1. \notag
  \end{align}
  If $s < r$, then $M_s$ is a submatrix of $M_r$ obtained by deleting
  a set of rows and corresponding columns, and therefore $\nu (M_s)
  \le \nu (M_r)$, yielding (1).
  
  (2) The first and third equalities hold by \eqref{E:stronger}.
  Since $J_rT_r^* = \tilde T_rJ_r$ and hence
  $$
  I_r - \tilde T_r\tilde T_r^* = I_r - J_rT_r^*T_rJ_r = J_r(I_r -
  T_r^*T_r)J_r,
  $$
  the second equality also holds.
  
  (3) By part (2), $\nu (I_n - \tilde T_n \tilde T_n^*) = \nu (I_n -
  T_nT_n^*) = \kappa$.  By the proof of (1), if $1 \le r \le n/2$,
  then
  $$
  \nu \begin{pmatrix} I_r-T_r T_r^* &Q_r \\ Q_r^*&I_r-\tilde
    T_r\tilde T_r^*
\end{pmatrix} = \nu (M_{2r}) 
\le \nu (M_n) = \nu (I_n - T_nT_n^*) = \kappa,
$$
and this proves (3).
\end{proof}

\begin{corollary}\label{C:equivalence}
  Let $a_0,a_1,a_2,\dots$ be complex numbers and define
  $c_0=1,c_1,c_2,\dots$ by \eqref{E:connection}.  If one of the three
  nondecreasing sequences
  \begin{equation}\notag
    \left\{
      \nu (I_r - T_r T_r^*)
    \right\}_1^\infty ,
\quad
    \left\{
      \nu (I_r - \tilde T_r \tilde T_r^*)
    \right\}_1^\infty ,
\quad
\left\{
  \nu \begin{pmatrix}
  I_r-T_rT_r^* &Q_r \\[3pt] Q_r^*&I_r-\tilde T_r\tilde T_r^*
\end{pmatrix}
\right\}_1^\infty 
  \end{equation}
  has constant value $\kappa$ from some point on, then all do.  If one
  of the three nondecreasing sequences
  \begin{equation}\notag
    \left\{
      \pi (I_r - T_r T_r^*)
    \right\}_1^\infty ,
\quad
    \left\{
      \pi (I_r - \tilde T_r \tilde T_r^*)
    \right\}_1^\infty ,
\quad
\left\{
  \pi \begin{pmatrix}
  I_r-T_rT_r^* &Q_r \\[3pt] Q_r^*&I_r-\tilde T_r\tilde T_r^*
\end{pmatrix}
\right\}_1^\infty 
  \end{equation}
  has constant value $\kappa$ from some point on, then all do.
\end{corollary}

\begin{proof}
  This follows on expressing all of the quantities in terms of the
  sequences $\{ \nu ( M_r) \}_{r=1}^\infty$ and $\{ \pi (M_r)
  \}_{r=1}^\infty$.  For example, for the negative eigenvalues, if one
  of the quantities has constant value $\kappa$ from some point on,
  then $\nu (M_r) = \kappa$ for all sufficiently large $r$, and all
  have constant value $\kappa$ from some point on.
\end{proof}

We next recall a result from \cite{CG2} on convergence of power
series. We give a complete proof not only to keep the note
selfcontained but also to show, as was done in Section 2, the role of
realization theory: the coefficients of the power series are
represented as Taylor coefficients of the transfer function of a
colligation, which is holomorphic in a neighborhood of the origin.

\begin{theorem}\label{T:CG0}
\label{T:redundant}
Let $a_0,a_1,a_2,\dots$ be complex numbers such that the matrices
$I_j-T_jT_j^*$ have $\kappa $ negative eigenvalues for all
sufficiently large $j$. Then the power series $S(z) =
\sum_{j=0}^\infty a_j z^j$ converges in some disk $|z|<\delta$ where
$\delta >0$.
\end{theorem}

\begin{proof}  Let $\eF = \boldC$ be the complex
  numbers viewed as a Hilbert space in the Euclidean metric.  Define
  $c_0,c_1,c_2,\dots$ by \eqref{E:connection}.  Then by
  \eqref{E:stronger}, the matrices \eqref{E:tmp1} have $\kappa$
  negative eigenvalues for all sufficiently large~$r$, that is, the
  sequence $c_0,c_1,c_2,\dots$ belongs to $\eP_\kappa$.  As in
  Iokhvidov and Kre{\u\i}n \cite[pp.~312-314]{MR21:6543}, construct a
  Naimark dilation for $c_0,c_1,c_2,\dots$; that is, we construct a
  Pontryagin space $\eK$ that contains $\eF$ isometrically as a
  regular subspace, and a unitary operator $U\in\lk$ such that
  $$
  c_j = P_\eF U^j\vert_\eF,\qquad j=0,1,2,\dots,
  $$
  where $P_\eF$ is the projection on $\eK$ with range~$\eF$.  Since
  $\eF$ is a regular subspace of $\eK$, we can write $ \eK = \eH
  \oplus \eF $ where $\eH$ is a regular subspace of $\eK$.  Let
\begin{equation}\notag
U = \begin{pmatrix}
A & B \\
C & D
\end{pmatrix}
\end{equation}
relative to this decomposition.  We show that
\begin{equation}\label{E:coefficients}
a_0 = D\qquad \text{and}\qquad a_m = CA^{m-1}B, \quad m \ge 1.
\end{equation}
The cases $m=0,1$ are immediate.  We prove the formula for $a_m$
assuming it is known for $a_0,\dots,a_{m-1}$.  By
\eqref{E:connection},
\begin{equation}\notag
  \label{E:connection2}
  c_{m+1} = c_0a_m + c_1a_{m-1} + \cdots + c_m a_0,
\end{equation}
so it is the same thing to show that
\begin{equation}
  \label{E:connection3}
  c_{m+1} = c_0 CA^{m-1}B + c_1 C A^{m-2}B + \cdots + c_{m-1}CB + c_mD.
\end{equation}
Put
\begin{equation}\notag
  \label{E:Uj}
   \begin{pmatrix}
     A & B \\
     C & D
\end{pmatrix}^j = \begin{pmatrix}
  A_j & B_j \\
  C_j & D_j
\end{pmatrix},
\qquad j \ge 0.
\end{equation}
Then
\begin{equation}\notag
  \label{E:Umplus1}
  \begin{pmatrix}
    A_{m+1} & B_{m+1} \\
    C_{m+1} & D_{m+1}
\end{pmatrix}
= \begin{pmatrix}
A_m & B_m \\
C_m & D_m
\end{pmatrix} \begin{pmatrix}
  A & B \\
  C & D
\end{pmatrix}
= \begin{pmatrix}
A_mA + B_mC & A_mB+B_mD \\
C_mA+D_mC & C_mB+D_mD
\end{pmatrix}.
\end{equation}
Since $D_j = P_\eF U^j \vert_\eF = c_j$ for all $j\ge 0$, $ c_{m+1} =
C_mB + c_mD.  $ This allows us to bring \eqref{E:connection3} to the
form
\begin{equation}
  \label{E:reduced}
  C_mB = c_0 CA^{m-1}B + c_1 C A^{m-2}B + \cdots + c_{m-1}CB.
\end{equation}
Dropping the factor $B$ on the right in each term, we easily verify
\eqref{E:reduced} by induction: the formula is evident for $m=1$, and
the inductive step follows from the identity $C_{m+1} = C_mA + c_mC$.
This completes the proof of \eqref{E:coefficients}.  The identity
\eqref{E:coefficients} implies that $|a_j| \le K \rho^j$ for some
positive constants $K$ and $\rho$, and therefore the power series
$\sum_{j=0}^\infty a_j z^j$ converges in a neighborhood of the origin.
\end{proof}

We can now relate Problem \ref{P:CF} to an indefinite form of the
trigonometric moment problem.

Let $\eP_\kappa$ ($\eP_{\nu, \pi}$) be the set of all sequences $\{
c_j \}_{j=0}^\infty$ with $\bar c_0 = c_0$ such that the matrix $M_r$
has $\kappa$ negative ($\nu$ negative and $\pi$ positive) eigenvalues
for all sufficiently large $r$.

\begin{problem}[Trigonometric moment problem]\label{P:TMP} 
  Let $c_0,c_1,\dots,c_{n-1}$ be $n$ complex numbers with $\bar c_0 =
  c_0$.  Determine for which nonnegative integers $\kappa$ there is a
  sequence $\{ c_p \}_{p=0}^\infty$ in $\eP_\kappa$ that extends the
  given numbers.  Determine for which nonnegative integers $\nu$ and
  $\pi$ there is a sequence $\{ c_p \}_{p=0}^\infty$ in $\eP_{\nu,
    \pi}$ that extends the given numbers.

\end{problem}

This problem is an indefinite form of the trigonometric moment problem
and it was considered by Iokhvidov and Kre{\u\i}n \cite[\S
19]{MR21:6543}.  In the classical case, this concerns the Fourier
coefficients, or moments, $c_j = \int e^{-ijt} \; d\mu(t)$, $j=0,\pm
1,\pm 2,\dots$, of a nonnegative measure $\mu$ on $[0,2\pi)$.  In this
case, the matrix $( c_{i-j} )_{i,j=0}^n$ is nonnegative for every $n
\ge 0$, since
\begin{equation}\notag
    \sum_{j,k=0}^n c_{k-j} \lambda_k \bar \lambda_j =
  \int_{[0,2\pi)}
  \sum_{j,k=0}^n \lambda_k \bar \lambda_j e^{-i(k-j)t} \; d\mu(t) 
= \int_{[0,2\pi)} \bigg| \sum_{j,k=0}^n
  \lambda_k e^{-ikt}\bigg|^2 \; d\mu(t) \ge 0
\end{equation}
for arbitrary numbers $\lambda_0,\dots,\lambda_n$.  When $\mu$ is a
probability measure, $c_0 = 1$.  The classical trigonometric moment
problem is to extend given numbers $c_0,c_1,\dots,c_{n-1}$ with $\bar
c_0 = c_0$ to such a moment sequence.  In the indefinite extension, we
still speak of the ``trigonometric moment problem,'' but the
underlying function theory is not the same.

We can show now that Problem \ref{P:CF} and Problem \ref{P:TMP} are
equivalent.

\begin{theorem}[Equivalence of Problems \ref{P:CF} and \ref{P:TMP}]
  \label{T:equiv}
  Assume that the numbers $a_0,\dots,a_{n-1}$ and $c_0 =
  1,c_1,\dots,c_{n-1},c_n$ are connected as in \eqref{E:connection}.
  Then Problem~{\ref{P:CF}} is solvable with the data
  $a_0,\dots,a_{n-1}$ if and only if Problem~{\ref{P:TMP}} is solvable
  with the data $c_0,\dots,c_{n-1},c_n$.
\end{theorem}

\begin{proof}
  Suppose that Problem~\ref{P:CF} with the data $a_0,\dots,a_{n-1}$
  has a solution in $\Skappa$. Let
  $$S(z) = \sum_{j=0}^\infty a_j z^j$$
  be the Taylor expansion of this
  solution.  By the necessary conditions for Problem~\ref{P:CF}
  discussed above, $I_j - T_jT_j^*$ has $\kappa$ negative eigenvalues
  for all sufficiently large $j$.  Define $c_{n+1},c_{n+2},\dots$ so
  that
  $$
  c_j = c_0a_{j-1}+c_1a_{j-2}+\cdots c_{j-1}a_0
  $$
  for all $j=1,2,\dots$.  Then \eqref{E:stronger} implies that
  $M_j$ has $\kappa$ negative eigenvalues for all $j=0,1,2,\dots$.
  Therefore $c_0,c_1,c_2,\dots$ is a solution to Problem~\ref{P:TMP}
  with the data $c_0,\dots,c_{n-1},c_n$.
  
  Conversely, assume that Problem~\ref{P:TMP} is solvable with the
  data $c_0,\dots,c_{n-1},c_n$, that is, the numbers can be extended
  to a sequence $c_0,c_1,c_2,\dots$ in $\eP_\kappa$.  Then the
  matrices \eqref{E:tmp1} have $\kappa$ negative eigenvalues for all
  sufficiently large $r$.  Reversing the process above, we obtain a
  sequence $a_0,a_1,a_2,\dots$ that extends $a_0,\dots,a_{n-1}$ such
  that the matrices $I_j - T_jT_j^*$ have $\kappa$ negative
  eigenvalues for all sufficiently large $j$.  By
  Theorem~\ref{T:redundant}, the series $S(z) = \sum_{j=0}^\infty a_j
  z^j$ converges in some disk $|z|<\delta$ where $\delta >0$, and by a
  theorem of Kre{\u\i}n and Langer in \cite[Theorem 6.3]{MR57:1173},
  the function $S(z)$ so defined belongs to $\Skappa$.  Thus
  Problem~\ref{P:CF} is solvable with the data $a_0,\dots,a_{n-1}$.
  The argument for the classes $\eP_{\nu,\pi}$ and ${\bf S}_{\nu,
    \pi}$ is similar.
\end{proof}

We use a series of propositions from \cite{Iokhvidov}.  The matrices
$M_0,M_1,M_2,\dots$ that appear in the list below are Hermitian
matrices of the form \eqref{E:tmp1} defined for appropriate numbers
$c_0 = \bar c_0, c_1, c_2, \dots$, and $n$ is any positive integer.
Recall that for any Hermitian matrix $A$ we write $\pi(A)$ and
$\nu(A)$ for the numbers of positive and negative eigenvalues of $A$
counting multiplicity.  The {\bf signature} of $A$ is $\sigma(A) =
\pi(A) - \nu(A)$.  Write $|A|$ for the determinant of $A$ and $\rho(A)
= \pi(A) + \nu(A)$ for the rank of $A$.

\begin{enumerate}
\item [$1^\circ$)] The difference $\rho(M_n)-\rho(M_{n-1})$ is either
  $0$, $1$, or $2$.  \smallskip
  
\item [$2^\circ$)] If $\rho(M_n)-\rho(M_{n-1})=0$, then $\pi(M_n)=
  \pi(M_{n-1})$ and $\nu(M_n)= \nu(M_{n-1})$.\smallskip
  
\item [$3^\circ$)] If $\rho(M_n)-\rho(M_{n-1})=1$, then either
  $\pi(M_n)= \pi(M_{n-1})+1 $ and $\nu(M_n)= \nu(M_{n-1})$, or
  $\pi(M_n)= \pi(M_{n-1})$ and $\nu(M_n)= \nu(M_{n-1})+1$ .\smallskip
  
\item [$4^\circ$)] If $\rho(M_n)-\rho(M_{n-1})=2$, then $\pi(M_n)=
  \pi(M_{n-1})+1$ and $\nu(M_n)= \nu(M_{n-1})+1$.  \smallskip
  
\item [$5^\circ$)] If $|M_{n-1}|\neq 0$, then there are infinitely
  many $c_n$ such that $\rho(M_n)=\rho(M_{n-1})$.  \smallskip
  
\item [$6^\circ$)] If $|M_{n-1}|= 0$ and $|M_{\rho(M_{n-1})-1}|\neq
  0$, then there is a unique $c_n$ such that $\rho(M_n)=\rho(M_{n-1})$
  .  \smallskip
  
\item [$7^\circ$)] The assumptions in $6^\circ$$)$ imply that there is
  a unique extension $(c_j)_{j=0}^{\infty}$ of $(c_j)_{j=0}^{n-1}$
  such that $\rho(M_j)=\rho(M_{n-1})$, $j\geq n$.\smallskip
  
\item [$8^\circ$)] There exists a $c_n$ with $\rho(M_n)=\rho(M_{n-1})$
  if and only if $|M_{\rho(M_{n-1})-1}|\neq 0$.  \smallskip
  
\item [$9^\circ$)] If $|M_{r-1}|\neq 0$ and $|M_{n-1}|=\cdots
  =|M_r|=0$ for some $0\leq r< \rho(M_{n-1})$ ($|M_{-1}|=1$ by
  definition), then $\rho(M_n)=\rho(M_{n-1})+2$.  \smallskip
  
\item [$10^\circ$)] If $|M_{n-1}|\neq 0$, then for each $k=1,2,\dots,$
  there are infinitely many $c_n$, $\dots$, $c_{n+k-1}$ such that
  $\nu(M_{n+k-1})=\nu(M_{n-1})+k$ and $|M_{n+k-1}|\neq 0$.\smallskip
  
\item [$11^\circ$)] If $|M_{n-1}|\neq 0$, then for each $\ell =1,2,
  \dots, $ there are infinitely many $c_n$, $\dots$, $c_{n+\ell -1}$
  such that $\pi(M_{n+\ell -1})=\pi(M_{n-1})+\ell $ and $|M_{n+\ell
    -1}|\neq 0$.\smallskip
  
\item [$12^\circ$)] $\sigma(M_{n-1})=\sum_{j=0}^{n-1} {\rm sign}\,
  (|M_{j-1}||M_j|)$, where, by definition, $|M_{-1}|=1$ and ${\rm
    sign}\,0=0$.  \smallskip
  
\item [$13^\circ$)] If $\rho(M_j)$ is a constant $\rho$ for all
  sufficiently large $j$, then $|M_{\rho -1}|\neq 0$.
\end{enumerate}

\begin{proof}[Proofs]  All of the citations below are from
  \cite{Iokhvidov}.
  
  $1^\circ$) Corollary on p.~34.
  
  $2^\circ$) Theorem 6.2, p.  36.
  
  $3^\circ$) Theorem 6.3, p.~36.
  
  $4^\circ$) Theorem 6.1, p.~35.
  
  $5^\circ$) Theorem 13.1, p.~97, and Remark 1, p.~98.
  
  $6^\circ$) Theorem 13.2, p.~100, and Remark 1, p.~102.
  
  $7^\circ$) Corollary on p.~101 and Remark 1, p.~102.
  
  $8^\circ$) The ``if'' part follows from $5^\circ)$ and $6^\circ$),
  the ``only if'' part from Theorem 15.3, p.~119.
  
  $9^\circ$) Proposition $3^\circ$, p.~121.
  
  $10^\circ$) and $11^\circ$) It is enough to prove these statements
  for $k=1$ in $10^\circ$) and $\ell=1$ in $11^\circ$).  To do this,
  we use the proof of Theorem 13.1, p.~97, and Remark 1, p.~98, to
  construct infinitely many extensions with $|M_n| > 0$ and infinitely
  many extensions with $|M_n| < 0$ (treat the subcases $|M_{n-2}| \neq
  0$ and $|M_{n-2}| = 0$ separately using the argument on p.~99).
  Then $10^\circ$) and $11^\circ$) follow from $3^\circ$).
  
  $12^\circ$) Theorem 16.1, p.~129.
  
  $13^\circ$) Theorem 15.4, p.~119.
\end{proof}

Our solution of Problem II is presented in Theorem~\ref{T:trig}.  The
first parts of the statements (a), (c), and (f) can be found in
Iokhvidov's book as Excercise~8 on pp.~133--134; in the interest of
completeness we prove these statements as well.  It is clear that a
given sequence $(c_j)_{j=0}^{n-1}$ does not have any extension
$(c_j)_{j=0}^{\infty}$ in $\eP_{\nu}$ if $\nu < \nu(M_{n-1})$, and
there is no extension in $\eP_{\nu, \pi}$ if either $\nu <
\nu(M_{n-1})$ or $\pi < \pi(M_{n-1})$, because by
$1^\circ$)--$4^\circ$), $\nu(M_j)$ and $\pi(M_j)$ are nondecreasing
functions of $j$.  If an extension $(c_j)_{j=0}^{\infty}$ belongs to
the class $\eP_{\nu}$ then it is possible that $\rho(M_j)$ and hence
also $\pi(M_j)$ tends to $\infty$ as $j \rightarrow \infty$.  Such an
extension does not belong to any of the classes $\eP_{\nu, \pi}$.
According to $13^\circ$) a necessary condition for
$(c_j)_{j=0}^{\infty}$ to belong to $\eP_{\nu, \pi}$ is that
$|M_{\nu+\pi-1}|\neq 0$.

\begin{theorem}\label{T:trig}
  Let $c_0 = \bar c_0, c_1, \dots, c_{n-1}$ be given numbers, and
  define $M_0,\dots,M_{n-1}$ as in~\eqref{E:tmp1}.

\noindent 
Assume $|M_{n-1}|\neq 0$.
  \begin{enumerate}
  \item[(a)] There exist infinitely many extensions in
    $\eP_{\nu(M_{n-1})}$, even infinitely many extensions in the
    smaller set $\eP_{\nu(M_{n-1}), \pi(M_{n-1})}$.  \smallskip
    
  \item[(b)] There exist infinitely many extensions in
    $\eP_{\nu(M_{n-1})+\nu, \pi(M_{n-1}) + \pi}$ for all $\nu \ge 0$
    and $\pi \ge 0$.
  \end{enumerate}
\noindent  Assume $|M_{n-1}|=0 $ and $|M_{\rho(M_{n-1})-1}|\neq 0$.
\begin{enumerate}
\item[(c)] There is a unique extension in $\eP_{\nu(M_{n-1})}$; it
  belongs to $\eP_{\nu(M_{n-1}), \pi(M_{n-1})}$.  \smallskip
  
\item[(d)] There are no extensions in $\eP_{\nu}$ for $\nu(M_{n-1})<
  \nu < \nu(M_{n-1})+{\rm dim}\,{\rm ker}\,M_{n-1};$ there are no
  extensions in $\eP_{\nu, \pi}$ if
  $$\nu(M_{n-1})< \nu < \nu(M_{n-1})+{\rm dim}\,{\rm ker}\,M_{n-1}$$
  or if
  $$\pi(M_{n-1})< \pi < \pi(M_{n-1})+{\rm dim}\,{\rm ker}\,M_{n-1}.$$
\item[(e)] There are infinitely many extensions in $\eP_{\nu, \pi}$
  for all pairs $(\nu,\pi)$ with
  $$\nu \geq \nu(M_{n-1})+{\rm dim}\,{\rm ker}\,M_{n-1} \quad
  \text{and} \quad \pi \geq \pi(M_{n-1})+{\rm dim}\,{\rm
    ker}\,M_{n-1}.$$
\end{enumerate}
\noindent Assume $|M_{n-1}|=0 $ and $|M_{\rho(M_{n-1})-1}|= 0$.
\begin{enumerate}
\item[(f)] There are no extensions in $\eP_{\nu(M_{n-1})}$.
  \smallskip
  
\item[(g)] There are no extensions in $\eP_{\nu}$ if $\nu <
  \nu(M_{n-1})+{\rm dim}\,{\rm ker}\,M_{n-1};$ there are no extensions
  in $\eP_{\nu, \pi}$ if
  $$
  \nu < \nu(M_{n-1})+{\rm dim}\,{\rm ker}\,M_{n-1}$$
  or if
  $$\pi < \pi(M_{n-1})+{\rm dim}\,{\rm ker}\,M_{n-1}.$$
\item[(h)] There are infinitely many extensions in $\eP_{\nu,\pi}$ for
  every pair $(\nu, \pi)$ with
  $$\nu \geq \nu(M_{n-1})+{\rm dim}\,{\rm ker}\,M_{n-1} \quad
  \text{and} \quad \pi \geq \pi(M_{n-1})+{\rm dim}\,{\rm
    ker}\,M_{n-1}.$$
\end{enumerate}
\end{theorem}

\begin{proof}   For any extension 
  of the given sequence by numbers $c_n,c_{n+1},\dots$, we assume that
  $M_n,M_{n+1},\dots$ are defined as in \eqref{E:tmp1}.
  
  (a) According to $5^\circ$) there are infinitely many $c_n$ such
  that $\rho(M_n)=\rho(M_{n-1})=n$. For such $M_n$ we have $|M_n|=0$
  and $|M_{\rho(M_n)-1}|\neq 0$. Hence by $7^\circ$) there is an
  extension $(c_j)_{j=0}^{\infty}$ of $(c_j)_{j=0}^{n-1}$ such that
  $\rho(M_j)=\rho(M_{n-1})$ for all $j \geq n-1$. Statement $2^\circ$)
  implies that
  $$\nu(M_j)=\nu(M_{n-1}) \quad \text{and} \quad
  \pi(M_j)=\pi(M_{n-1}), \qquad j\geq n-1,$$
  and hence
  $(c_j)_{j=0}^{\infty}$ belongs to $\eP_{\nu(M_{n-1}),
    \pi(M_{n-1})}$.
  
  (b) By $10^\circ$) there are infinitely many numbers $c_n$ such that
  $\nu(M_{n})=\nu(M_{n-1})+1$ and $|M_n|\neq 0$.  Therefore
  $\rho(M_n)=\rho(M_{n-1})+1$ and by $3^\circ$),
  $\pi(M_n)=\pi(M_{n-1})$.  After $\nu $ steps, we obtain numbers
  $c_n,\dots,c_{n+\nu -1}$ such that
  $$|M_{n+\nu -1}|\neq 0, \qquad \nu (M_{n+\nu -1})=\nu, \quad
  \text{and} \quad \pi(M_{n+\nu -1})=\pi(M_{n-1}).
  $$
  Using the same argument with $11^\circ$) instead of $10^\circ$),
  we obtain numbers $c_{n+\nu },\dots,c_{n+\nu +\pi -1}$ (each of
  which can be chosen in infinitely many ways) such that
  $$|M_{n+\nu +\pi -1}|\neq 0, \qquad \nu(M_{n+\nu +\pi -1})=\nu,
  \quad \text{and} \quad \pi(M_{n+\nu +\pi -1})=\pi.
  $$
  Now (b) follows from~(a).
  
  (c) According to $7^\circ$) there exists a unique extension
  $(c_j)_{j=0}^{\infty }$ of $(c_j)_{j=0}^{n-1}$ such that
  $$
  \rho(M_j)=\rho(M_{n-1}), \qquad j\geq n.
  $$
  It follows from $2^\circ$) that also $\pi(M_j)=\pi(M_{n-1})$ and
  $\nu(M_j)=\nu(M_{n-1})$ for $j\geq n$.  Therefore there exists a
  unique extension of $(c_j)_{j=0}^{n-1}$ in the class
  $\eP_{\nu(M_{n-1})}$ and this extension belongs to
  $\eP_{\nu(M_{n-1}), \pi(M_{n-1})}$ (for the uniqueness part, note
  that by $3^\circ$) the equality $\nu(M_n)=\nu(M_{n-1})$ can only
  hold in the present situation when $\rho(M_n)=\rho(M_{n-1})$).
  
  (d) and (e).  By hypothesis
  \begin{equation}\label{E:alpha}
    |M_{n-1}| = 0 \qquad \text{and} \qquad |M_{\rho(M_{n-1})-1}| \neq 0.
  \end{equation}
  The unique extension described in part (c) of the theorem cannot
  meet any of the conditions in parts (d) and (e); since for this
  extension $\rho(M_{n-1}) = \rho(M_{n}) = \rho(M_{n+1}) = \cdots$, in
  parts (d) and (e) we need only consider extensions such that
  \begin{equation}\notag
    \rho(M_{n-1}) = \cdots = \rho(M_{n+k-1})  < \rho(M_{n+k})
  \end{equation}
  for some $k \ge 0$.  In this situation \eqref{E:alpha} holds with
  $n$ replaced by $n+k$, and therefore we may restrict attention to
  extensions satisfying
  \begin{equation}\label{E:beta}
    \rho(M_{n-1}) < \rho(M_{n}).
  \end{equation}
  By $6^\circ)$, \eqref{E:beta} holds for all but one choice of $c_n$;
  in what follows, we assume that $c_n$ is chosen so that
  \eqref{E:beta} is satisfied.  The question then is if the sequence
  $(c_j)_{j=0}^n$ can be further extended to an infinite sequence
  $(c_j)_{j=0}^\infty$ as required in (d) and (e).

  \noindent Case (i): $\rho(M_n) = n+1$.
  
  Since $\rho(M_{n-1}) < n$ by \eqref{E:alpha}, by $1^\circ)$ we must
  have $\rho(M_{n-1}) = n-1$.  Thus ${\rm dim}\,{\rm ker}\,M_{n-1} =
  1$, and hence part (d) holds vacuously.  Part (e) also holds in this
  case.  For by statement $4^\circ$), $\nu(M_n)=\nu(M_{n-1})+1$ and
  $\pi(M_n)=\pi(M_{n-1})+1$ and since $M_n$ is invertible, part (e)
  follows from (a).
  
  \goodbreak
  \noindent Case (ii): $\rho(M_n) < n+1$.
  
  Then with $r=\rho(M_{n-1})$, in view of \eqref{E:alpha} and
  \eqref{E:beta},
  $$|M_{r-1}|=|M_{\rho(M_{n-1})-1}|\neq 0, \quad |M_r|=\cdots
  =|M_{n-1}|=|M_n|=0.$$
  Consider any extension of $(c_j)_{j=0}^n$ by a
  number $c_{n+1}$.  By \eqref{E:beta}, $r < \rho(M_{n+1})$.  Applying
  $9^\circ)$ with $n$ replaced by $n+1$, we obtain
  $$
  \rho(M_{n+1}) = \rho(M_{n})+2,
  $$
  and by $4^\circ)$,
  $$
  \nu(M_{n+1})=\nu(M_n)+1 \qquad \text{and} \qquad
  \pi(M_{n+1})=\pi(M_n)+1.
  $$
  If $\rho(M_{n+1}) < n+2$, we can repeat this argument.  We
  continue in this way for $k = 1,2,\dots$ and extend
  $(c_j)_{j=0}^{n}$ with any numbers $c_{n+1}, \dots, c_{n+k}$,
  $k=1,2, \dots$; by $9^\circ$) and $4^\circ$), we have $r <
  \rho(M_{n+k})$,
  \begin{align}
    \rho(M_{n+k}) &= \rho(M_n)+2k, \notag \\
    \nu(M_{n+k})  &= \nu(M_n)+k,   \notag \\
    \pi(M_{n+k}) &= \pi(M_n)+k, \notag
  \end{align}
  and
  $$|M_{r-1}|\neq 0, \quad |M_r|=\cdots =|M_{n-1}|=|M_n|=\cdots
  =|M_{n+k}|=0,$$
  provided $\rho(M_{n+k})=\rho(M_n)+2k <n+k+1$. If
  equality holds, that is,
  $$k=k_0:= n-\rho(M_n)+1,$$
  then $M_{n+k_0}$ is invertible and the
  process stops. Hence if such an extension of $(c_j)_{j=0}^n$ can be
  continued to a sequence in a class $\eP_{\nu, \pi}$, then
  necessarily
  $$\nu \geq \nu_0:=\nu(M_{n+k_0})=\nu(M_n)+k_0=
  \nu(M_n)+n-\rho(M_n)+1 ,$$
  $$\pi \geq \pi_0:=\pi(M_{n+k_0})=\pi(M_n)+k_0=
  \pi(M_n)+n-\rho(M_n)+1 ,$$
  and according to (a) and (b) each of the
  classes $\eP_{\nu}$ and $\eP_{\nu,\pi}$ contains infinitely many
  extensions. Thus the first part of (d) and (e) will follow once we
  show that
  $$\nu_0 = \nu(M_{n-1})+{\rm dim}\,{\rm ker}\,M_{n-1} \quad
  \text{and} \quad \pi_0 = \pi(M_{n-1})+{\rm dim}\,{\rm ker}\,M_{n-1}.
  $$
  Since $|M_{n-1}|=0$, $12^\circ$) implies that
  $$\sigma(M_n)- \sigma(M_{n-1})={\rm sign}\,|M_{n-1}||M_n|=0$$
  and
  since $\rho(M_n)> \rho(M_{n-1})$, we therefore have $\rho(M_n)=
  \rho(M_{n-1})+2$, and by $4^\circ$), $\nu(M_n)=\nu(M_{n-1})+1$ and
  $\pi(M_n) =\pi(M_{n-1})+1$.  This implies that $\nu_0=
  \nu(M_{n-1})+{\rm dim}\,{\rm ker}\,M_{n-1}$ and also that $\pi_0$
  has the desired value.
  
  From the first part of (d) it follows that there are no extensions
  in $\eP_{\nu, \pi}$ if
  $$\nu(M_{n-1}) < \nu < \nu(M_{n-1})+{\rm dim}\,{\rm ker}\,M_{n-1},$$
  whatever the value of $\pi$. By considering the sequence
  $(-c_j)_{j=0}^{n-1}$ and its extensions $(-c_j)_{j=0}^{\infty}$ and
  applying the results just proved (together with
  $\nu(-M_j)=\pi(M_j)$) we find that there are no extensions in
  $\eP_{\nu, \pi}$ if
  $$\pi(M_{n-1}) < \pi < \pi(M_{n-1})+{\rm dim}\,{\rm ker}\,M_{n-1}$$
  whatever the value of $\nu$.
  
  \goodbreak (f) is part of (g).
  
  (g) and (h). By $8^\circ)$, \eqref{E:beta} holds for any choice of
  $c_n$.  This allows us to proceed by an argument which is similar to
  the proof of (d) and (e) above; in case (ii) there, the exact value
  of $r$ is unimportant in order to obtain the conclusion.
\end{proof}

We can now deal with Problem~\ref{P:CF}.  According to Theorem
\ref{T:equiv}, we must apply the previous result to the case where
$c_0=1$, $n$ is replaced by $n+1$, $|M_n|=|I_n-T_nT_n^*|$, ${\rm
  dim}\,{\rm ker}\,M_n={\rm dim}\,{\rm ker}\,(I_n-T_nT_n^*)$, and
$$\rho(I_n-T_nT_n^*)=\rho(M_n)-1, \ \pi(I_n-T_nT_n^*)=\pi(M_n)-1, \ 
\nu(I_n-T_nT_n^*)=\nu(M_n).
$$
Note that $\eP_{\nu, \pi}$ corresponds to the class ${\bf S}_{\nu,
  \pi'}$ with $\pi'=\pi-1$.  We obtain the following solution for the
Carath\'eodory-Fej\'er problem.

\begin{theorem}\label{T:CF} 
  Let $a_0,\dots,a_{n-1}$ be given numbers, and define
  $T_1,\dots,T_{n}$ as in $\eqref{E:toeplitz1}$.

\noindent
Assume $|I_n-T_nT^*_n|\neq 0$.
  \begin{enumerate}
  \item[(a$'$)] There exist infinitely many solutions of Problem
    \ref{P:CF} in ${\bf S}_{\nu(I_n-T_nT^*_n)}$, even in the smaller
    set ${\bf S}_{\nu(I_n-T_nT^*_n), \pi(I_n-T_nT^*_n)}$.  \smallskip
    
  \item[(b$'$)] There exist infinitely many solutions in ${\bf
      S}_{\nu, \pi}$ for all pairs $(\nu, \pi)$ with $\nu \geq
    \nu(I_n-T_nT^*_n)$ and $\pi \geq \pi(I_n-T_nT^*_n)$.
  \end{enumerate}
\noindent Assume $|I_n-T_nT^*_n|=0 $ and $|I_{\rho }-T_{\rho }T^*_{\rho }|\neq
0$, where $\rho =\rho (I_n-T_nT^*_n)$.
  \begin{enumerate}
  \item[(c$'$)] There is a unique solution in ${\bf
      S}_{\nu(I_n-T_nT^*_n)}$; it belongs to ${\bf
      S}_{\nu(I_n-T_nT^*_n), \pi(I_n-T_nT^*_n)}$.  \smallskip
    
  \item[(d$'$)] There are no solutions in ${\bf S}_{\nu}$ for
    $\nu(I_n-T_nT^*_n)< \nu < \nu(I_n-T_nT^*_n)+{\rm dim}\,{\rm
      ker}\,I_n-T_nT^*_n;$ there are no solutions in ${\bf S}_{\nu,
      \pi}$ if
    $$\nu(I_n-T_nT^*_n)< \nu < \nu(I_n-T_nT^*_n)+{\rm dim}\,{\rm
      ker}\,(I_n-T_nT^*_n)$$
    or if
    $$\pi(I_n-T_nT^*_n) < \pi < \pi(I_n-T_nT^*_n)+{\rm dim}\,{\rm
      ker}\,(I_n-T_nT^*_n).$$
  \item[(e$'$)] There are infinitely many solutions in ${\bf S}_{\nu,
      \pi}$ for all pairs $(\nu ,\pi)$ with
    $$\nu \geq \nu(I_n-T_nT^*_n)+{\rm dim}\,{\rm ker}\,(I_n-T_nT^*_n)
    \quad \text{and} \quad \pi \geq \pi(I_n-T_nT^*_n)+{\rm dim}\,{\rm
      ker}\,(I_n-T_nT^*_n).$$
  \end{enumerate}
\noindent Assume $|I_n-T_nT^*_n|=0 $ and 
$|I_{\rho }-T_{\rho }T^*_{\rho }|= 0$.
\begin{enumerate}
\item[(f$'$)] There are no solutions in ${\bf S}_{\nu(I_n-T_nT^*_n)}$.
  \smallskip
  
\item[(g$'$)] There are no solutions in ${\bf S}_{\nu}$ if $\nu <
  \nu(I_n-T_nT^*_n)+{\rm dim}\,{\rm ker}\,(I_n-T_nT^*_n);$ there are
  no solutions in ${\bf S}_{\nu, \pi}$ if
  $$
  \nu < \nu(I_n-T_nT^*_n)+{\rm dim}\,{\rm ker}\,(I_n-T_nT^*_n)$$
  or
  if
  $$\pi < \pi(I_n-T_nT^*_n)+{\rm dim}\,{\rm ker}\,(I_n-T_nT^*_n).$$
\item[(h$'$)] There are infinitely many solutions in ${\bf
    S}_{\nu,\pi}$ for every pair $(\nu, \pi)$ with
  $$\nu \geq \nu(I_n-T_nT^*_n)+{\rm dim}\,{\rm ker}\,(I_n-T_nT^*_n)
  \quad \text{and} \quad \pi \geq \pi(I_n-T_nT^*_n)+{\rm dim}\,{\rm
    ker}\,(I_n-T_nT^*_n).$$
\end{enumerate}
\end{theorem}

\goodbreak We mention a consequence of the solution of
Problem~\ref{P:CF} for the case $\nu =\nu (I_n-T_nT^*_n)$.

\begin{corollary}
  Let $a_0,a_1,\dots,a_{n-1}$ be numbers such that $I_n - T_nT_n^*$
  has $\nu$ negative eigenvalues.
  \begin{enumerate}
  \item[(1)] If $I_n-T_nT_n^*$ is invertible, Problem~\ref{P:CF} has
    infinitely many solutions in ${\bf S}_{\nu }$.  \smallskip
    
  \item[(2)] If $I_n-T_nT_n^*$ is singular and $\rho (I_{n-1} -
    T_{n-1}T_{n-1}^*) = \rho (I_{n} - T_{n}T_{n}^*)$,
    Problem~\ref{P:CF} has a unique solution in ${\bf S}_{\nu }$.
    \smallskip
    
  \item[(3)] If $I_n-T_nT_n^*$ is singular and $\rho (I_{n-1} -
    T_{n-1}T_{n-1}^*) < \rho (I_{n} - T_{n}T_{n}^*)$,
    Problem~\ref{P:CF} has no solution in ${\bf S}_{\nu }$.
  \end{enumerate}
\end{corollary}

The results in \cite{CG1} and \cite{CG2} give a solution to the
existence and uniqueness problems for the matrix versions of both the
trigonometric moment problem and Carath\'eodory-Fej\'er coefficients
problem in the indefinite case, but the question of obtaining a matrix
analogue of (d$'$)-(h$'$) in Theorem \ref{T:CF} is open.

\section{A remark on holomorphic kernels}\label{S:kernels}

The result below is used in Section~\ref{S:coefficients} and is well
known in particular cases.  The general result is presumably also
known, but we do not know a reference.  For the convenience of the
reader, we sketch a proof.

Let $K(w,z)= \sum_{m,n=0}^\infty C_{mn} z^m \bar w^n $ be a
holomorphic Hermitian kernel defined for $|w|<R$ and $|z|<R$, with
values in $\lf$ for some Kre{\u\i}n space $\eF$.  For any nonnegative
integer $r$, we may alternatively view the matrix $(C_{mn})_{m,n=0}^r$
as a selfadjoint operator on $\eF^{r+1} = \eF \oplus \cdots \oplus
\eF$, where there are $r+1$ summands on the right side, or as a kernel
on a finite set.  The number of negative eigenvalues of
$(C_{mn})_{m,n=0}^r$ as an operator and the number of negative squares
of $(C_{mn})_{m,n=0}^r$ as a kernel coincide.

\begin{theorem}\label{negsquares}
  Let $\kappa$ be a nonnegative integer.  Then $\sqminus K = \kappa$
  if and only if
  $$\nu (C_{mn})_{m,n=0}^r \le \kappa
  $$
  for all nonnegative integers $r$ and equality holds for all
  sufficiently large $r$.
\end{theorem}

We can formulate this result in another way.  Let ${\mathbb N}_0$ be
the set of nonnegative integers.  Define a kernel $C$ on ${\mathbb
  N}_0 \times {\mathbb N}_0$ by
$$
C(m,n) = C_{mn}, \qquad m,n \in {\mathbb N}_0.
$$
Then $\sqminus K = \sqminus C$.  The theory of Kolmogorov
decompositions \cite{cgkolmogorov} gives a natural approach to this
result, but we base our argument on similar notions for reproducing
kernel Pontryagin spaces.

\begin{proof}
  Since a holomorphic Hermitian kernel has the same number of negative
  squares on subregions \cite[Theorem 1.1.4]{ADRSbook}, by a change of
  scale we may assume that $R > 1$.  We may also assume without loss
  of generality that $\eF$ is a Hilbert space.  Let $H^2_{\eF}$ be the
  Hardy class of $\eF$-valued functions on the unit disk $\boldD$.
  
  Assume that $\sqminus K = \kappa$.  By a method of Alpay
  \cite{MR92m:46038}, we define a bounded selfadjoint operator $P$ on
  $H^2_{\eF}$ such that
  $$
  P \colon (1 - \bar w z)^{-1}f \to K(w,z)f, \qquad w \in \boldD, f
  \in \eF,
  $$
  and
  $$
  P \colon z^nf \to A_n(z)f, \qquad f \in\eF, \; n = 0,1,2,\dots.
  $$
  where $K(w,z) = \sum_{n=0}^\infty A_n(z) \bar w^n$, that is,
  $A_n(z) = \sum_{m=0}^\infty C_{mn} z^m$ for all $n = 0,1,2,\dots$.
  For another account of the construction of $P$, see \cite[Theorem
  8.4]{JRTopics}.  By the spectral theorem, we can write $P = P_+ +
  P_0 + P_-$, where $P_{\pm}$ and $P_0$ are selfadjoint operators
  corresponding to the spectral subspaces $\eH _+$, $\eH_-$, and
  $\eH_0 = \ker P$ for the sets $(0,\infty)$, $(-\infty,0)$, and
  $\{0\}$.  Since $\sqminus K = \kappa$, $\dim \eH_- = \kappa$.  Let
  $\eK_0$ be $H^2_\eF/\ker P$.  Write $[h] = h + \ker P$ for the coset
  determined by an element $h$ of $H^2_\eF$.  Define a nondegenerate
  inner product on $\eK_0$ by
  $$
  {\ip{[h]}{[k]}}_{\eK_0} = {\ip{Ph}{k}}_{H^2_\eF}, \qquad h,k \in
  H^2_\eF.
  $$
  Using \cite[Theorem 2.5, p.~20]{IKL}, complete $\eK_0$ to a
  Pontryagin space $\eK$ having negative index~$\kappa$.  The cosets
  determined by the polynomials are dense in $H^2_{\eF}/\ker P$ by
  \cite[statement (i) on p.~20]{IKL}, and therefore $\{ [z^nf] \colon
  f \in \eF, n = 0,1,2,\dots \}$ is a total set in $\eK$. By
  construction,
  $$
  {\ip{[z^mf_1]}{[z^nf_2]}}_{\eK} = {\ip{C_{mn}f_1}{f_2}}_{\eF},
  \qquad f_1,f_2 \in \eF,\; m,n = 0,1,2,\dots .
  $$
  Hence by \cite[Lemma 1.1.1]{ADRSbook}, the matrix
  $(C_{mn})_{m,n=0}^r$ has at most $\kappa$ negative eigenvalues for
  all $r=0,1,2,\dots$ and one such matrix has exactly $\kappa$
  negative eigenvalues.  Since the number of negative eigenvalues of
  $(C_{mn})_{m,n=0}^r$ is a nondecreasing function of $r$, this number
  is $\kappa$ for all sufficiently large $r$.
  
  Conversely, assume that the matrix $(C_{mn})_{m,n=0}^r$ has at most
  $\kappa$ negative eigenvalues for all $r=0,1,2,\dots$ and exactly
  $\kappa$ negative eigenvalues for all sufficiently large~$r$.  By
  what we showed above, if we can only show that $\sqminus K \le
  \kappa$, it will follow that $\sqminus K = \kappa$.  Let ${\mathbb
    N}_0$ be the set of nonnegative integers, and define a kernel $C$
  on ${\mathbb N}_0 \times {\mathbb N}_0$ by
  $$
  C(m,n) = C_{mn}, \qquad m,n \in {\mathbb N}_0.$$
  Our hypotheses
  imply that $\sqminus C = \kappa$.  By \cite[Theorem
  1.1.3]{ADRSbook}, there is a unique Pontryagin space $\eH_C$ of
  functions $h = \{ h_n \}_{n=0}^\infty$ on ${\mathbb N}_0$ with
  reproducing kernel $C$.  This means that for each $m \in {\mathbb
    N}_0$ and $f \in\eF$, the sequence $C(m,\cdot)f = \{ C_{mn}f
  \}_{n=0}^\infty$ belongs to $\eH_C$, and for any element $h = \{ h_n
  \}_{n=0}^\infty$ of $\eH_C$,
  $$
  {\ip{\{ h_n \}_{n=0}^\infty}{\{ C_{mn}f \}_{n=0}^\infty}}_{\eH_C}
  = {\ip{h_m}{f}}_{\eF}.
  $$
  By \cite[Theorem 1.1.2]{ADRSbook}, we can represent the kernel
  $C$ in the form
  $$
  C_{mn} = A_n^*A_m, \qquad m,n \in {\mathbb N}_0,
  $$
  where for each $k \in {\mathbb N}_0$, $A_k^*$ is the evaluation
  mapping on $\eH_C$ to $\eF$: $ A_k^* \left( \{ h_n \}_{n=0}^\infty
  \right) = h_k.  $ By the Cauchy representation, the operators
  $C_{mn}$ are uniformly bounded, and therefore for $w$ and $z$ in a
  suitable neighborhood of the origin,
  $$
  K(w,z) = \sum_{m,n=0}^\infty A_n^*A_m z^m \bar w^n = A(w)^*A(z),
  $$
  where $A(z) = \sum_{m=0}^\infty A_m z^m$.  The values of $A(z)$
  lie in the Pontryagin space $\eH_C$, which has negative
  index~$\kappa$.  The restriction of $K(w,z)$ to a suitable
  neighborhood of the origin thus has at most $\kappa$ negative
  squares, and since the number of negative squares is independent of
  the domain (see \cite[Theorem 1.1.4]{ADRSbook}), $\sqminus K \le
  \kappa$.  As noted above, this implies that $\sqminus K = \kappa$.
\end{proof}

\providecommand{\bysame}{\leavevmode\hbox to3em{\hrulefill}\thinspace}

\noindent 2000 Mathematics Subject Classification:
Primary 47A57, 42A70, 30E05

\hskip7.23cm Secondary 47B50, 47A48, 47A06

\end{document}